\documentclass[12pt]{amsart}

\usepackage{amsmath}
\usepackage{amssymb}
\usepackage{amsthm}
\usepackage{amscd}
\usepackage{epsfig}
\usepackage{epic}
\usepackage{eepic}

\newcommand     {\comment}[1]   {}
\newcommand{\mute}[2] {}
\newcommand     {\printname}[1] {}
\newcommand{\labell}[1] {\label{#1}\printname{#1}}

\let\mathbbm\mathbb

\newtheorem{theorem}{Theorem}
\newtheorem{Theorem}[theorem]{Theorem}
\newtheorem{proposition}[theorem]{Proposition}
\newtheorem{Lemma}[theorem]{Lemma}
\newtheorem*{Lemma*}{Lemma}

\newtheorem{Corollary}[theorem]{Corollary}
\newtheorem{Proposition}[theorem]{Proposition}

\theoremstyle{remark}

\newtheorem{Remark}[theorem]{Remark}
\newtheorem*{Remark*}{Remark}

\newtheorem{Example}[theorem]{Example}
\newtheorem*{Example*}{Example}

\theoremstyle{definition}

\newtheorem{Definition}[theorem]{Definition}
\newtheorem*{Definition*}{Definition}
\newtheorem*{Notation*}{Notation}

\newcommand{\function}[3]{\ensuremath{#1\colon #2\to #3}}
\newcommand{\R}{\ensuremath{\mathbbm{R}}}

\newcommand{\Z}{\ensuremath{\mathbbm{Z}}}
\newcommand{\N}{\ensuremath{\mathbbm{N}}}

\newcommand{\GL}[2]{\ensuremath{{\rm GL}_{#1}(#2)}}

\newcommand{\eval}[1]{\ensuremath{\lvert_{#1}}}

\DeclareMathOperator{\interior}{interior}

\newcommand{\fF}{\ensuremath{\mathfrak F}}

\newcommand{\calU}{\ensuremath{\mathcal U}}

\def \calF {\mathcal{F}}
\def \calR {\mathcal{R}}
\def \inv  {^{-1}}
\def \ol   {\overline}
\def \eps  {\epsilon}
\def \ssminus  {\smallsetminus}

\DeclareMathOperator{\image}{image}
\def \sss {\scriptstyle}

\def \tU {\widetilde{U}}
\def \tu {\widetilde{u}}
\def \pibar {\ol{\pi}}

\def\scal<#1,#2>{\langle #1\,,#2\rangle}

\title{Orbifolds as diffeologies}

\author{Patrick Iglesias}
\address{CNRS, France, and The Hebrew University of Jerusalem, Israel}
\email{piz@math.huji.ac.il}

\author{Yael Karshon}
\address{The University of Toronto}
\email{karshon@math.toronto.edu}

\author{Moshe Zadka}
\address{The Hebrew University of Jerusalem}
\email{moshez@divmod.com}

\thanks{This research was partially supported by the Israel Science
Foundation founded by the Academy of Sciences and Humanities, by the
National Center for Scientific Research (CNRS, France),
and by the National Science and Engineering Research Council 
of Canada (NSERC)}

\begin{document}

\begin{abstract}
We consider orbifolds as diffeological spaces.
This gives rise to a natural notion of differentiable maps
between orbifolds, making them into a subcategory of diffeology.
We prove that the diffeological approach to orbifolds
is equivalent to Satake's notion of a V-manifold
and to Haefliger's notion of an orbifold.
This follows from a lemma: a diffeomorphism (in the diffeological sense)
of finite linear quotients lifts to an equivariant diffeomorphism.
\end{abstract}


\maketitle


\section{Introduction}
\labell{sec:intro}

An orbifold is a space that is locally modeled on quotients
of $\R^n$ by finite linear group actions.
The precise definition of its global structure is more complicated.
The situation in the literature is somewhat problematic:
different authors give definitions that are a priori different
from each other, and the relations between them are not made
clear.

We propose to approach orbifolds through the notion of a \emph{diffeology}.
A diffeology on a set $X$ specifies, for each whole number $m$ 
and each open subset $U$ of $\R^m$, which functions from $U$ to $X$
are differentiable.  A \emph{diffeological space} is a set equipped
with a diffeology.  
A quotient of a diffeological space is again a diffeological space; 
a map to the quotient is declared to be differentiable if and only if 
it locally lifts.  We define a diffeological orbifold to be 
a diffeological space that is locally diffeomorphic at each point 
to the quotient of $\R^n$ by a finite linear group action.
For details, see Section \ref{sec:diffeology}.

In Section~\ref{sec:lifting} we show that a diffeomorphism
(in the diffeological sense) between finite linear quotients
of $\R^n$ lifts to an equivariant diffeomorphism.
See Lemma~\ref{crucial}. 
This fails for arbitrary differentiable maps;
see Examples~\ref{example1} and~\ref{example2}.

To illustrate that the diffeology retains rich information,
in Section~\ref{sec:structure groups} we show that the
stabilizers of a finite group action are a diffeological
invariant of the quotient.

Orbifolds were originally introduced by Satake under the name
``V-manifolds".  
Satake's definition is subtle and is often quoted imprecisely.  
We recall this definition in Section~\ref{sec:Vmanifolds}.
A problem with the definition is that it does not lead to a
satisfactory notion of $C^\infty$ maps.
(In \cite[page 469]{satake2}, Satake writes this footnote:
``The notion of $C^\infty$-map thus defined is inconvenient in the point
that a composite of two $C^\infty$-maps defined in a different choice
of defining families is not always a $C^\infty$ map.'')
A related problem with Satake's definition 
is that it is not a priori clear 
whether a V-manifold structure is \emph{local}, that is, whether 
a V-manifold structure can be determined by its restrictions 
to the elements of an open covering.
In diffeology, the notion of a differentiable map is completely 
natural, and, with it, diffeological orbifolds become a category.
Also, the axioms of diffeology immediately imply that a diffeological
orbifold structure is local.

In order to not increase the confusion in the literature,
we verify that, under mild restrictions 
(Satake insisted that the space be Hausdorff and the finite groups
act without reflections), our definition of a diffeological orbifold
is equivalent to Satake's notion of a V-manifold.
See Theorem \ref{main theorem}.
The proof relies on the technical results of Section~\ref{sec:lifting}.
As an application, we prove, in Section \ref{sec:locality},
that a V-manifold structure is local.
In section~\ref{sec:Haefliger} we similarly prove 
that a diffeological orbifold is the same thing as an orbifold structure
in the sense of Haefliger~\cite{haefliger}.

The modern literature contains ``higher" approaches to orbifolds,
for example, through stacks or Lie groupoids.
See \cite{DM,M:orbifold,M,BX}.
In these approaches, the finite group actions can be non-effective;
the resulting structure is global, not local.
A concrete relationship between the Lie groupoid approach
to orbifolds and the diffeology approach to orbifolds
is analyzed in an forthcoming paper by Y.~Karshon and M.~Zoghi.

We expect our diffeological approach to orbifolds
to suffice for applications that do not involve
global aspects of non-effective orbifolds.
It is extremely easy to work with diffeology.  Thus,
when applicable, diffeology removes the language barrier
to passing from manifolds to orbifolds.

\subsubsection*{Acknowledgement}
P.~I.\ writes: I'm grateful to the
Hebrew University of Jerusalem for its hospitality that has lead us
to this collaboration.
Y.~K.\ writes: I am grateful to Sue Tolman for convincing me, years ago,
that I don't understand orbifolds.  In her 1993 thesis \cite{tolman:thesis}
she defined a cohomology theory for orbifolds and showed how to compute it
for moduli spaces of complex curves.

\section{Diffeological orbifolds} 
\labell{sec:diffeology}

The notions of diffeology and diffeological spaces
are due to Chen \cite{Chen1} (under the name ``differentiable spaces") 
and to Souriau \cite{Souriau1,Souriau2}.
For details, see \cite{Iglesias1,Iglesias:web}.

A \emph{parametrization} on a set $X$ is a map \function{p}{U}{X},
where $U$ is an open subset of $\R^n$ for some $n$.

\begin{Definition}[Diffeology and diffeological spaces]
\labell{def:diffeology} A \emph{diffeology} on a set $X$ is a set of
parametrizations on $X$, whose elements are called \emph{plots}, such
that the following axioms hold:

\begin{description}
\item [Covering axiom]
	The constant parametrizations are plots.
\item[Locality axiom] Let $p \colon U \to X$ be a parametrization.
	If every point in $U$ is contained in an open subset $V \subset U$
	such that $p|_V \colon V \to X$ is a plot, then $p$ is a plot.
\item[Smooth compatibility axiom]
	Let $\function{p}{U}{X}$ be a plot, $V$ an open subset
	of $\R^m$ for some $m$, and $F \colon V \to U$ a $C^\infty$
         smooth map.  Then $p \circ F$ is a plot.
\end{description}
A \emph{diffeological space} is a set $X$ equipped with a diffeology.
\end{Definition}

\begin{Definition}[Differentiable maps] \labell{def:differentiable}
A map \function{f}{X}{Y} between diffeological spaces is called
\emph{differentiable} if the composition of $f$ with any plot of $X$
is a plot of $Y$.
A map \function{f}{X}{Y} between diffeological spaces is called
a \emph{diffeomorphism} if $f$ is a bijection and both $f$ and $f^{-1}$
are differentiable.
\end{Definition}

Diffeological spaces with differentiable maps form a category.

\smallskip

Like the notion of a topology,
the notion of a diffeology is very general.  For instance, the
set of \emph{all} maps $U \to X$ satisfies the three diffeology axioms;
it defines the \emph{coarse diffeology} on $X$.  The set of
all locally constant maps $U \to X$ defines the \emph{discrete
diffeology}
on $X$.  Other diffeologies are described below.
In each of the following definitions, we leave it to the reader
to verify that the three diffeology axioms hold.

\begin{Example*}[Smooth manifolds]
On a manifold $M$, the set of all $C^\infty$ smooth maps
$p \colon U \to M$ is a diffeology.
A map between manifolds is differentiable
as a map of diffeological spaces
if and only if it is smooth in the usual sense.
\end{Example*}

\begin{Remark}
Any plot $U \to X$ of a diffeological space $X$ is differentiable,
where the open subset $U \subset \R^n$ is taken with its manifold 
diffeology.  
This is the content of the smooth compatibility axiom.
\end{Remark}

\begin{Definition}[Quotient diffeology]\labell{def:quotientD}
Let $X$ be a diffeological space, let $\sim$ be an equivalence relation
on $X$, and let $\pi \colon X \to Y:=X/\!\sim$ be the quotient map.
The \emph{quotient diffeology} on $Y$
is the diffeology in which \function{p}{U}{Y} is a plot
if and only if each point in $U$ has a neighborhood $V \subset U$
and a plot \function{\Tilde{p}}{V}{X}
such that $p\eval{V} = \pi \circ \Tilde{p}$.
\end{Definition}

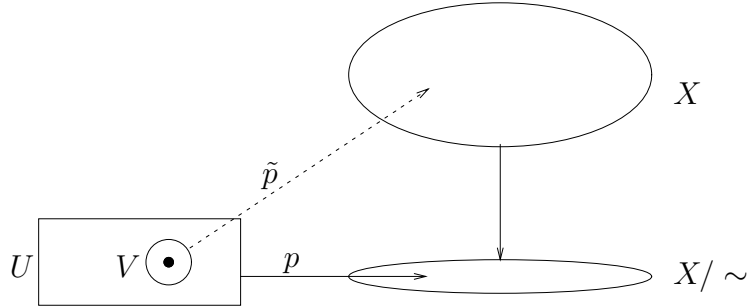
\begin{figure}[ht]
\begin{center}
\setlength{\unitlength}{0.0005in}
\begingroup\makeatletter\ifx\SetFigFont\undefined%
\gdef\SetFigFont#1#2#3#4#5{%
  \reset@font\fontsize{#1}{#2pt}%
  \fontfamily{#3}\fontseries{#4}\fontshape{#5}%
  \selectfont}%
\fi\endgroup%
{\renewcommand{\dashlinestretch}{30}
\begin{picture}(7361,3034)(0,-10)
\put(1650,462){\blacken\ellipse{100}{100}}
\put(1650,462){\ellipse{100}{100}}
\put(1650,462){\ellipse{474}{474}}
\put(5100,312){\ellipse{3150}{350}}
\put(5100,2412){\ellipse{3150}{1500}}
\path(300,912)(2400,912)(2400,12)
	(300,12)(300,912)
\path(5100,1687)(5100,487)
\path(5070.000,607.000)(5100.000,487.000)(5130.000,607.000)
\path(2400,312)(4325,312)
\path(4205.000,282.000)(4325.000,312.000)(4205.000,342.000)
\dashline{45.000}(1875,612)(4350,2262)
\path(4266.795,2170.474)(4350.000,2262.000)(4233.513,2220.397)
\put(0,312){\makebox(0,0)[lb]{\smash{{{\SetFigFont{12}{14.4}{\rmdefault}{\mddefault}{\updefault}$U$}}}}}
\put(1100,312){\makebox(0,0)[lb]{\smash{{{\SetFigFont{12}{14.4}{\rmdefault}{\mddefault}{\updefault}$V$}}}}}
\put(6900,237){\makebox(0,0)[lb]{\smash{{{\SetFigFont{12}{14.4}{\rmdefault}{\mddefault}{\updefault}$X/\sim$}}}}}
\put(6900,2112){\makebox(0,0)[lb]{\smash{{{\SetFigFont{12}{14.4}{\rmdefault}{\mddefault}{\updefault}$X$}}}}}
\put(2850,412){\makebox(0,0)[lb]{\smash{{{\SetFigFont{12}{14.4}{\rmdefault}{\mddefault}{\updefault}$p$}}}}}
\put(2625,1287){\makebox(0,0)[lb]{\smash{{{\SetFigFont{12}{14.4}{\rmdefault}{\mddefault}{\updefault}$\tilde{p}$}}}}}
\end{picture}
}
\caption{Quotient diffeology}
\label{fig:quotient}
\end{center}
\end{figure}

\begin{Example*}
Let a finite group $\Gamma$ act on $\R^n$ by linear transformations.
Then the quotient $\R^n / \Gamma$ is naturally a diffeological space.
\end{Example*}

\begin{Definition}[Local diffeomorphisms] \labell{local diffeomorphism}
A diffeological space $X$ is \emph{locally diffeomorphic}
to a diffeological space $Y$ at a point $x \in X$ if and only if
there exists a subset $A$ of $X$, containing $x$, and there exists
a one-to-one function $f \colon A \to Y$ such that
\begin{enumerate}
\item
for any plot $p \colon U \to X$,
the composition $f \circ p$ is a plot of $Y$;
\item
for any plot $q \colon V \to Y$, the composition
$f\inv \circ q$ is a plot of $X$.
\end{enumerate}
\end{Definition}

\begin{Remark*}
The condition that the maps $f \circ p$ and $f\inv \circ q$ be plots
contains the requirement
that the domains of these maps, $p\inv(A)$ and $q\inv(f(A))$, be open.

Also see Remark \ref{loc diffeo II}.
\end{Remark*}

\begin{Example*}
An $n$ dimensional manifold can be re-defined as a diffeological space
which is locally diffeomorphic to $\R^n$ at each point.
Notice that being a manifold is a property of the pre-existing
diffeology; there is no need for any additional structure.
\end{Example*}

Taking a similar approach to the notion of orbifolds, we make the
following definition:

\begin{Definition}[Diffeological orbifolds] \labell{def:Dorbifold}
A \emph{diffeological orbifold} is a diffeological space
which is locally diffeomorphic at each point to a quotient
$\R^n/\Gamma$, for some $n$, where $\Gamma$ is a finite group
acting linearly on $\R^n$.
\end{Definition}

\begin{Remark*}
Diffeological orbifolds, with differentiable maps,
form a subcategory of the category of diffeological spaces.
\end{Remark*}

\begin{Definition}[D-topology] \labell{def:D topology}
The \emph{D-topology} on a diffeological space $X$
is the topology in which a subset $W \subset X$ is open
if and only if $p\inv W$ is open for every plot $p \colon U \to X$.
This is the finest topology on $X$ for which all plots are continuous.
A differentiable map between diffeological spaces is always continuous
with respect to their D-topologies. That is, we have a forgetful functor
from diffeological spaces to topological spaces. 
\end{Definition}

\begin{Example*} \labell{manifold topology}
On a smooth manifold, considered as a diffeological space,
the D-topology coincides with the standard topology.
On $\R^n / \Gamma$, where $\Gamma \subset \GL{n}{\R}$
is a finite subgroup, the D-topology coincides with the
standard quotient topology. 
\end{Example*}

\begin{Definition}[Subset diffeology and diffeological subspaces]
Let $X$ be a diffeological space, $A \subset X$ a subset,
and \function{i}{A}{X} the inclusion map.
Declare $p \colon U \to A$ to be a plot
if and only if $i \circ p$ is a plot of $X$.
This defines the \emph{subset diffeology} on $A$.
With this diffeology, $A$ is called
a \emph{diffeological sub-space of $X$}.
\end{Definition}

In particular, any D-open subset of a diffeological space
is naturally a diffeological space.

\begin{Remark} \labell{loc diffeo II}
A diffeological space $X$ is locally diffeomorphic to a diffeological
space $Y$ at a point $x \in X$
if and only if there exist a $D$-open subset $A$ of $X$ containing $x$
and a $D$-open subset $B$ of $Y$, equipped with their subset diffeologies,
and a diffeomorphism $f \colon A \to B$.
\end{Remark}

If $\Gamma \subset \GL{n}{\R}$ is a finite subgroup,
the open subsets of $\R^n/\Gamma$ are exactly the quotients $U/\Gamma$
where $U \subset \R^n$ is a $\Gamma$-invariant open set.
The diffeology on $U / \Gamma$ can be considered 
either as the quotient diffeology of $U$ or, equivalently, as the 
subset diffeology of $\R^n / \Gamma$.  Hence,

\begin{Remark} \labell{part of orbifold}
Any D-open subset of a diffeological orbifold is
also a diffeological orbifold.
\end{Remark}


\begin{Example} \labell{finite stabs}
Let a finite group $G$ act on a manifold $M$.
Then the quotient $M/G$ is a diffeological orbifold.
More generally, if a Lie group $G$ acts on a manifold $M$ 
properly and with finite stabilizers, then the quotient $M/G$ 
is a diffeological orbifold.
\end{Example}

\begin{proof}
This is an immediate consequence of the slice theorem,
by which a neighborhood of each $G$-orbit is equivariantly diffeomorphic
to $G \times_H D$, where $H \subset G$ is the stabilizer of a point 
in the orbit, acting linearly on a disk $D$.  See \cite{koszul} 
for compact (in particular, finite) group actions and \cite{palais} 
for proper group actions.
\end{proof}

\begin{Lemma} \labell{D cover}
Let $X$ be a diffeological space and $\calU$ an open cover of $X$ 
with respect to the D-topology.  Then the diffeology of $X$ 
is uniquely determined by the diffeology of the sets $U \in \calU$.
\end{Lemma}

\begin{proof}
This follows from Definition \ref{def:D topology} and the locality
axiom of a diffeology. 
\end{proof}

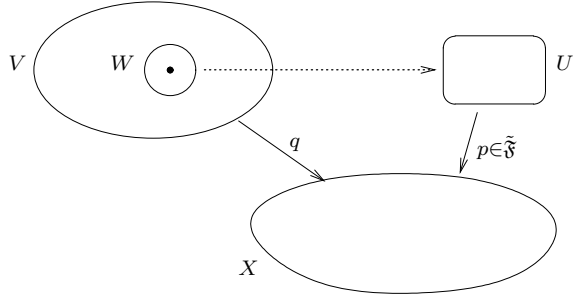
\begin{figure}[ht]
\begin{center}
\setlength{\unitlength}{0.0006in}
\begingroup\makeatletter\ifx\SetFigFont\undefined%
\gdef\SetFigFont#1#2#3#4#5{%
  \reset@font\fontsize{#1}{#2pt}%
  \fontfamily{#3}\fontseries{#4}\fontshape{#5}%
  \selectfont}%
\fi\endgroup%
{\renewcommand{\dashlinestretch}{30}
\begin{picture}(4830,2601)(0,-10)
\put(1425,1979){\ellipse{450}{450}}
\put(1275,1979){\ellipse{2100}{1200}}
\put(1425,1979){\blacken\ellipse{50}{50}}
\put(1425,1979){\ellipse{50}{50}}
\put(3930,1784){\arc{210}{1.5708}{3.1416}}
\put(3930,2174){\arc{210}{3.1416}{4.7124}}
\put(4620,2174){\arc{210}{4.7124}{6.2832}}
\put(4620,1784){\arc{210}{0}{1.5708}}
\path(3825,1784)(3825,2174)
\path(3930,2279)(4620,2279)
\path(4725,2174)(4725,1784)
\path(4620,1679)(3930,1679)
\put(4820,1979){\makebox(0,0)[lb]{\smash{{{\SetFigFont{12}{14.4}{\rmdefault}{\mddefault}{\updefault}$\sss U$}}}}}
\dottedline{45}(1725,1979)(3750,1979)
\path(3630.000,1949.000)(3750.000,1979.000)(3630.000,2009.000)
\path(2025,1529)(2775,1004)
\path(2659.488,1048.239)(2775.000,1004.000)(2693.896,1097.392)
\path(4125,1604)(3975,1079)
\path(3979.121,1202.625)(3975.000,1079.000)(4036.812,1186.141)
\path(2888,1029)(2917,1033)(2947,1037)
	(2978,1041)(3009,1044)(3042,1047)
	(3075,1050)(3109,1053)(3145,1056)
	(3181,1058)(3218,1060)(3255,1061)
	(3293,1063)(3332,1064)(3371,1065)
	(3410,1065)(3450,1065)(3490,1065)
	(3529,1065)(3568,1064)(3607,1063)
	(3645,1061)(3682,1060)(3719,1058)
	(3755,1056)(3791,1053)(3825,1050)
	(3858,1047)(3891,1044)(3922,1041)
	(3953,1037)(3983,1033)(4013,1029)
	(4048,1024)(4082,1018)(4116,1011)
	(4150,1004)(4183,997)(4216,989)
	(4249,981)(4281,972)(4313,963)
	(4344,953)(4375,943)(4404,933)
	(4433,922)(4461,911)(4488,900)
	(4513,889)(4537,878)(4560,866)
	(4582,855)(4603,844)(4622,833)
	(4640,822)(4656,811)(4672,800)
	(4686,790)(4700,779)(4715,766)
	(4729,754)(4742,741)(4755,728)
	(4766,714)(4776,700)(4785,686)
	(4794,671)(4801,656)(4807,641)
	(4812,626)(4815,610)(4817,594)
	(4818,578)(4818,562)(4817,547)
	(4815,531)(4811,516)(4807,500)
	(4801,485)(4795,470)(4788,454)
	(4780,440)(4771,425)(4762,410)
	(4751,394)(4739,378)(4725,362)
	(4710,345)(4693,328)(4675,311)
	(4655,294)(4633,276)(4610,259)
	(4585,242)(4559,225)(4531,209)
	(4502,193)(4472,178)(4440,164)
	(4407,150)(4373,137)(4338,124)
	(4302,113)(4264,102)(4225,91)
	(4196,85)(4166,78)(4136,72)
	(4104,65)(4070,59)(4036,54)
	(4001,49)(3964,44)(3926,39)
	(3888,34)(3848,30)(3807,27)
	(3766,23)(3724,20)(3681,18)
	(3638,16)(3595,14)(3551,13)
	(3508,12)(3465,12)(3422,12)
	(3380,12)(3338,13)(3297,15)
	(3257,17)(3217,19)(3179,22)
	(3142,25)(3105,29)(3070,33)
	(3035,38)(3002,43)(2969,48)
	(2938,54)(2900,62)(2863,70)
	(2827,79)(2792,89)(2758,100)
	(2724,111)(2690,123)(2657,136)
	(2625,150)(2593,164)(2563,178)
	(2533,194)(2504,209)(2476,225)
	(2449,241)(2423,258)(2399,274)
	(2375,291)(2354,307)(2333,323)
	(2314,339)(2296,355)(2279,370)
	(2264,385)(2250,400)(2236,414)
	(2224,428)(2213,441)(2200,457)
	(2189,472)(2178,488)(2169,503)
	(2160,518)(2153,533)(2146,548)
	(2141,564)(2137,579)(2134,594)
	(2132,609)(2131,624)(2132,638)
	(2134,653)(2137,667)(2141,681)
	(2146,694)(2153,707)(2160,720)
	(2169,732)(2178,744)(2189,756)
	(2200,768)(2213,779)(2225,790)
	(2238,800)(2253,811)(2268,822)
	(2285,833)(2303,844)(2323,855)
	(2344,866)(2366,878)(2390,889)
	(2415,900)(2441,911)(2469,922)
	(2497,933)(2527,943)(2557,953)
	(2588,963)(2619,972)(2651,981)
	(2684,989)(2717,997)(2750,1004)
	(2784,1011)(2818,1018)(2852,1024)(2888,1029)
\put(0,1979){\makebox(0,0)[lb]{\smash{{{\SetFigFont{12}{14.4}{\rmdefault}{\mddefault}{\updefault}$\sss V$}}}}}
\put(900,1979){\makebox(0,0)[lb]{\smash{{{\SetFigFont{12}{14.4}{\rmdefault}{\mddefault}{\updefault}$\sss W$}}}}}
\put(2025,179){\makebox(0,0)[lb]{\smash{{{\SetFigFont{12}{14.4}{\rmdefault}{\mddefault}{\updefault}$\sss X$}}}}}
\put(2475,1304){\makebox(0,0)[lb]{\smash{{{\SetFigFont{12}{14.4}{\rmdefault}{\mddefault}{\updefault}$\sss q$}}}}}
\put(4125,1229){\makebox(0,0)[lb]{\smash{{{\SetFigFont{12}{14.4}{\rmdefault}{\mddefault}{\updefault}$\sss p \in \Tilde{\fF}$}}}}}
\end{picture}
}
\caption{Generating family}
\label{fig:generating}
\end{center}
\end{figure}

\begin{Definition}
Let $\fF$ be a family of parametrizations on a set $X$.
Let $\Tilde{\fF}$ be the union of $\fF$ with the set of all
constant parametrizations on $X$.  Declare $q \colon V \to X$
to be a plot if and only if each point in $V$ has a neighborhood
$W \subset V$, a parametrization $p \colon U \to X$ in
$\Tilde{\fF}$,
and a smooth map $\function{F}{W}{U}$, such that
$q\eval{W} = p \circ F$.
This defines the \emph{diffeology generated by $\fF$}.
See Figure \ref{fig:generating}.
\end{Definition}

We end this section with a characterization of diffeological orbifolds:

\begin{Lemma} \labell{centered}
A diffeological space $X$ is a diffeological orbifold
if and only if for every $x \in X$ and open neighborhood $O$ of $x$
there exists a finite group $\Gamma$ acting on $\R^n$ linearly
and a map $\varphi \colon \R^n \to X$ whose image is contained in $O$,
such that $\varphi(0) = x$,
and such that $\varphi$ induces a diffeomorphism of $\R^n / \Gamma$
with a neighborhood of $x$ in $X$.
\end{Lemma}

Lemma \ref{centered} is a consequence of the slice theorem.  
See Definition \ref{def:Dorbifold} and the proof of 
Example \ref{finite stabs}.
For completeness, let us also give a more elementary proof of this
fact.

For any finite subgroup $\Gamma \subset \GL{n}{\R}$,
let $\| \cdot \|_\Gamma$ be the norm on $\R^n$ that is induced
from the inner product
$\sum_{\gamma \in \Gamma} \left< \gamma u , \gamma v \right>$
where $\left< \cdot, \cdot \right>$ is the standard inner product.
Then $\| \cdot \|_\Gamma$ is $\Gamma$-invariant. 

\begin{Lemma} \labell{centered for Rn}
Let $\Gamma$ be a finite group, acting linearly on $\R^n$.
Let $u \in \R^n$ be a point and let
$\Gamma_u = \{ \gamma \in \Gamma \ | \ \gamma \cdot u = u \}
   \subset \Gamma$ be its stabilizer. Let $U\subset \R^n$ be an open
$\Gamma$-invariant subset containing $u$. Then there exists a
$\Gamma_u$-equivariant open map
$$ \Tilde{f} \colon \R^n \to U $$
that sends the origin to $u$, is a diffeomorphism with its image,
and such that the induced map
$$ \R^n / \Gamma_u \to U / \Gamma $$
is a diffeomorphism with its image.
\end{Lemma}

\begin{figure}[ht]
\begin{center}
\setlength{\unitlength}{0.00083333in}
\begingroup\makeatletter\ifx\SetFigFont\undefined%
\gdef\SetFigFont#1#2#3#4#5{%
  \reset@font\fontsize{#1}{#2pt}%
  \fontfamily{#3}\fontseries{#4}\fontshape{#5}%
  \selectfont}%
\fi\endgroup%
{\renewcommand{\dashlinestretch}{30}
\begin{picture}(2327,2177)(0,-10)

\put(1855,596){\blacken\ellipse{50}{50}}
\put(1855,596){\ellipse{50}{50}}
\put(1855,596){\ellipse{900}{900}}
\put(1655,596){\makebox(0,0)[lb]{\smash{{{\SetFigFont{12}{14.4}{\rmdefault}{\mddefault}{\updefault}$\sss u$}}}}}
\put(1855,0){\makebox(0,0)[lb]{\smash{{{\SetFigFont{12}{14.4}{\rmdefault}{\mddefault}{\updefault}$\sss B(u,\eps)$}}}}}

\put(511,511){\blacken\ellipse{50}{50}}
\put(511,511){\ellipse{50}{50}}
\put(511,511){\ellipse{900}{900}}
\put(136,511){\makebox(0,0)[lb]{\smash{{{\SetFigFont{12}{14.4}{\rmdefault}{\mddefault}{\updefault}$\sss \gamma_2 \cdot u$}}}}}

\put(1111,1711){\blacken\ellipse{50}{50}}
\put(1111,1711){\ellipse{50}{50}}
\put(1111,1711){\ellipse{900}{900}}
\put(736,1711){\makebox(0,0)[lb]{\smash{{{\SetFigFont{12}{14.4}{\rmdefault}{\mddefault}{\updefault}$\sss \gamma_1 \cdot u$}}}}}

\end{picture}
}
\end{center}
\caption{Neighborhood of a $\Gamma$-orbit}
\label{fig-Gamma_u}
\end{figure}
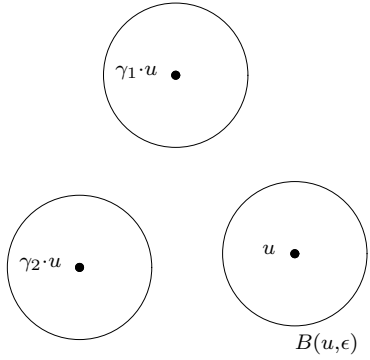

\begin{proof}
Let $\eps > 0$ be sufficiently small so that
$$ B(u,\eps) := \{ v \in \R^n \mid \| v - u \|_{\Gamma} < \eps \}$$
is contained in $U$ and so that $2\eps < \| \gamma \cdot u - u \|$
for all $\gamma \in \Gamma \ssminus \Gamma_u$.
Then the balls $B(u,\eps)$ and $\gamma \cdot B(u,\eps)$ 
are equal if $\gamma \in \Gamma_u$ and disjoint 
if $\gamma \in \Gamma \ssminus \Gamma_u$.
See Figure \ref{fig-Gamma_u},
where $\gamma_i$ represent the non-trivial cosets in $\Gamma/\Gamma_u$.
It follows that the inclusion map
$ B(u,\epsilon) \hookrightarrow \R^n $
induces a one to one map
\begin{equation} \labell{one2one}
 B(u,\epsilon)/\Gamma_u \to U/\Gamma .
\end{equation}
The vertical arrows in the diagram
$$
\begin{array}{ccc}
  B(u,\eps) & \hookrightarrow & U \\
\downarrow & & \downarrow \\
  B(u,\eps)/\Gamma_u & \to & U/\Gamma
\end{array}
$$
generate the quotient diffeologies;
it follows that the map \eqref{one2one} is a diffeomorphism with its image.
Let
$$ s \colon \R^n \to B(0,\epsilon)$$
be a $\Gamma_u$-equivariant diffeomorphism.  The map
$$ \Tilde{f}(x) = u + s(x) $$
has the desired properties.
\end{proof}

%
%
%
%

\begin{proof}[Proof of Lemma \ref{centered}]
The lemma follows immediately from Definition \ref{def:Dorbifold}
and Lemma \ref{centered for Rn}.
\end{proof}

\section{Lifting diffeomorphisms of finite quotients}
\labell{sec:lifting}

In this section we show that a diffeomorphism 
(in the diffeological sense)
of quotient spaces 
whose domain is $\R^n/\Gamma$, where $\Gamma \subset \GL{n}{\R}$
is a finite subgroup, always lifts to an equivariant diffeomorphism.

We note that this is not always true for differentiable maps
that are not diffeomorphisms:

\begin{Example}
Let $\Z/2\Z$ act on $\R^2$ by $(x,y) \mapsto \pm (x,y)$.
Consider a map from $\R^2$ to $\R^2/(\Z/2\Z)$ that is given by
$$(r \cos \theta , r \sin \theta) \mapsto
  [g(r) \cos(\theta/2) , g(r) \sin(\theta/2)],$$
where  $g$ is a function that vanishes near $r=0$
but does not vanish everywhere.  This map is well defined,
is differentiable, but does not lift to a 
differentiable (or even continuous)
map from $\R^2$ to $\R^2$.
\end{Example}

Also see examples~\ref{example1} and~\ref{example2}.

The following result is a crucial ingredient in our analysis:

\begin{Lemma} \labell{lem:identity}
Let $\Gamma \subset \GL{n}{\R}$ be a finite subgroup.
Let $U \subset \R^n$ be a connected open subset.  Let
$$ \function{h}{U}{\R^n} $$
be a continuously differentiable map that sends each point
to a point in the same $\Gamma$-orbit.
Then there exists a unique element $\gamma \in \Gamma$ such that
\begin{equation} \labell{h x gamma x}
    h(x) = \gamma x
\end{equation}
for all $x \in U$.
\end{Lemma}

\begin{proof}
For each $x \in U$ there exists $\gamma \in \Gamma$
such that \eqref{h x gamma x} holds.
We need to show that it is possible to choose the same $\gamma$
for all $x$.

Let $U_0$ denote the set of points in $U$ whose $\Gamma$-stabilizer
is trivial.  For $x \in U_0$, the element $\gamma \in \Gamma$ 
that satisfies \eqref{h x gamma x} is unique.
Because $h$ is continuous, this element is the same for all the points
in the same connected component of $U_0$.  For each connected component
$A$ of $U_0$, denote by $\gamma_A$ this element of $\Gamma$, so that
$$ h(x) = \gamma_A x \quad \text{ for all } \ x \in A.$$
Then, because this is a linear map,
$$ d_x h = \gamma_A \quad \text{ for all } \ x \in A,$$
and because $h$ is continuously differentiable,
$$ d_x h = \gamma_A \quad \text{ for all } \ x \in \text{closure}(A).$$
Because $U_0$ is obtained from $U$ by removing a finite union
of proper subspaces, namely, the fixed point sets of the elements
$\gamma \in \Gamma$ other than the identity, $U_0$ is dense in $U$.
Because $U$ is connected, it follows that all the $\gamma_A$'s 
are equal to each other.
\end{proof}

\begin{Remark}
In the above lemma, if $\Gamma$ contains no reflections,
it is enough to assume that $h$ is continuous:
if $\Gamma$ contains no reflections, the set $U_0$ is connected,
so there exists $\gamma$ such that \eqref{h x gamma x} holds
for all $x \in U_0$.  Because $U_0$ is dense and $h$ is continuous,
\eqref{h x gamma x} holds for all $x \in U$.
\end{Remark}

We also need the following topological fact:

\begin{Lemma} \labell{lem:U connected}
Let $\Gamma \subset \GL{n}{\R}$ be a finite subgroup.
Let $U \subset \R^n$ be a $\Gamma$-invariant subset that contains the 
origin.
If $U/\Gamma$ is connected, $U$ is connected.
\end{Lemma}

\begin{proof}
Suppose that $W \subset U$ is both open and closed in $U$.
After possibly switching between $W$ and $U \ssminus W$
we may assume that $0$ is in $W$.
Then each of the sets $\gamma \cdot W$, for $\gamma \in \Gamma$,
contains $0$ and is both open and closed in $U$.  The intersection
\begin{equation} \labell{intersection}
\bigcap\limits_{\gamma \in \Gamma} \gamma \cdot W
\end{equation}
is $\Gamma$-invariant, and its image in $U/\Gamma$
is both open and closed and is non-empty.  Because $U/\Gamma$
is connected, this image is all of $U/\Gamma$. Because the intersection
(\ref{intersection}) is $\Gamma$-invariant and its image is all of
$U/\Gamma$, it follows that this intersection is equal to $U$.
Thus, $W=U$.
\end{proof}

The following three lemmas show that a diffeomorphism of finite
quotients locally lifts to equivariant diffeomorphisms.

\begin{Notation*}
We will denote the image of the origin under the quotient map
$\R^n \to \R^n / \Gamma$ by the symbol $0$ and will still call it 
``the origin".
\end{Notation*}

\begin{Lemma}   
\labell{local liftings}
Let $\Gamma \subset \GL{n}{\R}$ and
$\Gamma' \subset \GL{n'}{\R}$ be finite subgroups.
Let $U \subset \R^n$ and $U' \subset \R^{n'}$
be invariant open subsets that contain the origin.  Let
$$ \function{f}{U/\Gamma}{U'/\Gamma'} $$
be a diffeomorphism of diffeological spaces such that $f(0) = 0$.  Let
$$ \Tilde{f} \colon U \to U' $$
be a smooth map that lifts $f$.  Then there exists a neighborhood
of the origin
on which $\Tilde{f}$ is a diffeomorphism with an open subset
of $U'$.  (In particular, it follows that $n = n'$.)
\end{Lemma}

\begin{proof}
By the definitions of the quotient diffeologies and of ``diffeomorphism"
(see Definitions \ref{def:differentiable} and \ref{def:quotientD}),
there exists an open $\Gamma'$-invariant subset $V' \subset U'$,
containing the origin, and a smooth map $\Tilde{g} \colon V' \to U$
that lifts the inverse map $f\inv \colon U'/\Gamma' \to U/\Gamma$.  
We may choose $V'$ to be connected, for instance,
by taking it to be a ball around the origin with respect to
a $\Gamma'$-invariant metric.
Let $V$ be the preimage of $f\inv(V'/\Gamma')$ under the quotient map
$U \to U/\Gamma$.
Then $V/\Gamma$ is connected, because it is the image of the connected
set $V'/\Gamma'$ under the continuous map $f\inv$.
By Lemma \ref{lem:U connected}, $V$ is also connected.
We now have connected invariant open neighborhoods $V$ and $V'$ 
of the origins in $\R^n$ and $\R^{n'}$ and smooth maps
$$ \Tilde{f} \colon V \to V' \quad \text{ and } \quad
     \Tilde{g} \colon V' \to V $$
such that $\Tilde{f}$ lifts $f$ and $\Tilde{g}$ lifts $f\inv$.
Without loss of generality we assume that $U=V$ and $U'=V'$.
We will show that $\Tilde{f}$ is then a diffeomorphism.

The composition $\Tilde{g} \circ \Tilde{f} \colon V \to V$
sends each $\Gamma$-orbit to itself.  By Lemma \ref{lem:identity}
there exists an element $\gamma$ of $\Gamma$ such that
$\Tilde{g} \circ \Tilde{f}(x) = \gamma x$ for all $x \in V$.
After replacing $\Tilde{f}$ by $\Tilde{f} \circ \gamma\inv$,
we may assume that
$$ \Tilde{g} \circ \Tilde{f}|_V \ \text{ is the identity map. }$$
Similarly, there exists an element $\gamma' \in \Gamma'$ such that
\begin{equation} \labell{with gamma}
   \Tilde{f} \circ \Tilde{g}(x) = \gamma' x \quad
   \text{ for all } x \in V'.
\end{equation}
Because $\Tilde{g} \circ \Tilde{f}$ is the identity,
we also have that
\begin{equation} \labell{without gamma}
   \Tilde{f} \circ \Tilde{g}(x) = x
   \quad \text{ for all $x \in V'$ in the image of $\Tilde{f}$ },
\end{equation}
so that
\begin{equation} \labell{gamma prime x}
   \gamma' x = x \quad
   \text{for all $x \in V'$ in the image of $\Tilde{f}$. }
\end{equation}
Because the image of $f$ is all of $V'/\Gamma'$,
the image of $\Tilde{f}$ contains an element $x$
whose stabilizer is trivial.
By this and \eqref{gamma prime x}, $\gamma' = 1$.
So $\Tilde{f}$ and $\Tilde{g}$ are inverses of each other.
In particular, $\Tilde{f}$ is a diffeomorphism.
\end{proof}

\begin{Lemma} 
\labell{lem:loc diffeo}
Let $\Gamma \subset \GL{n}{\R}$ and $\Gamma' \subset \GL{n'}{\R}$ 
be finite subgroups.  Let $U \subset \R^n$
and $U' \subset \R^{n'}$ be invariant open subsets.  Let
$$ f \colon U / \Gamma \to U' / \Gamma' $$
be a diffeomorphism of diffeological spaces.  Let
$$ \Tilde{f} \colon U \to U' $$
be a smooth map that lifts $f$.  
Then each point of $U$ has a neighborhood on which $\Tilde{f}$
is a diffeomorphism with an open subset of $U'$.

Consequently, $n=n'$, and for each $u \in U$ 
the linear map $L = d_u \Tilde{f}$ is invertible.
The conjugation map $\gamma \mapsto L \gamma L\inv$
carries the stabilizer subgroup $\Gamma_u$ of $u$
to the stabilizer subgroup $\Gamma'_{u'}$ of $u' = \Tilde{f}(u)$;
in particular, the stabilizer subgroups
$\Gamma_u$ and $\Gamma'_{u'}$ are conjugate in $\GL{n}{\R}$.
\end{Lemma}

\begin{proof}
Let $u\in U$ and $u'=\Tilde f(u)$.  By Lemma \ref{centered for Rn},
there exists a $\Gamma_u$-equivariant map $\Tilde g \colon \R^n\to U$
that sends the origin to $u$, is a diffeomorphism with an open
subset of $U$, and such that the induced map 
$g \colon \R^n/\Gamma_u \to U/\Gamma$ 
is a diffeomorphism with an open subset of $U/\Gamma$.
Similarly, there exists a $\Gamma'_{u'}$-equivariant map 
$\Tilde{g}' \colon \R^{n'}\to U'$ that sends the origin to $u'$, 
is a diffeomorphism with an open subset of $U'$, 
and such that the induced map $g' \colon \R^{n'}/\Gamma'_{u'} \to U'/\Gamma'$
is a diffeomorphism with an open subset of $U'/\Gamma'$.
The composition $(\Tilde g')^{-1}\circ \Tilde f \circ \Tilde g$
is a smooth map from a $\Gamma_u$-invariant open subset of $\R^n$
that contains the origin to a $\Gamma_{u'}$-invariant open subset
of $\R^{n'}$ that contains the origin.  It sends $0$ to $0$,
and it lifts the diffeomorphism $ (g')\inv \circ f \circ g$.
The first part of the lemma follows by applying Lemma \ref{local liftings}
to this composition.
\end{proof}

\begin{Lemma} \labell{lifting is eq diff}
Let $\Gamma \subset \GL{n}{\R}$ and $\Gamma' \subset \GL{n'}{\R}$
be finite subgroups. 
Let $U \subset \R^n$ be an invariant connected open subset
that contains the origin.  Let
$$\function{f}{U/\Gamma}{\R^{n'}/\Gamma'}$$
be a diffeomorphism of $U/\Gamma$ with an open subset of $\R^{n'}/\Gamma'$.
Let
$$\function{\Tilde{f}}{U}{\R^{n'}}$$
be a smooth map that lifts $f$. Then
\begin{enumerate}
\item
$\Tilde{f}$ is one to one.
\item
Let $U''$ denote the preimage of $f(U/\Gamma)$ under the projection map 
$\R^{n'} \to \R^{n'}/\Gamma'$; note that $U'' \subset \R^{n'}$ is open.
Then $\Tilde{f}$ is a diffeomorphism of $U$ with a connected component 
of $U''$.
\item
Let $L = d_0 \Tilde{f}$. 
Note that, by Lemma \ref{lem:loc diffeo}, 
$ h(\gamma) = L \gamma L\inv $
defines an isomorphism from $\Gamma$ to a subgroup of $\Gamma'$. 
Then $\Tilde{f}$ is $h$-equivariant.
\end{enumerate}
\end{Lemma}

\begin{proof}[Proof of (1)]
Because $\Tilde{f}$ lifts $f$ and $f$ is one to one, 
if $\Tilde{f}(u) = \Tilde{f}(v)$ then there exists $\gamma \in \Gamma$ 
such that $v = \gamma \cdot u$.  
We need to show, for each $\gamma \in \Gamma$, 
that $\tilde{f}(u) = \Tilde{f}(\gamma \cdot u)$
implies $\gamma \cdot u = u$.

Consider the set
$$ U_\gamma := \{ u \in U \ | \
     \Tilde{f}(\gamma \cdot u) = \Tilde{f}(u) \}.$$
Clearly, $U_\gamma$ is closed in $U$ and contains the fixed point set
$$ U^\gamma = \{ u \in U \ | \ \gamma \cdot u = u \} .$$
We need to show that $U_\gamma = U^\gamma$.

Let $v$ be any point in $U_\gamma$, and let
$\Gamma_v = \{ g \in \Gamma \ | \ g \cdot v = v \}$
be its stabilizer.
\begin{itemize}
\item
Suppose $\gamma \in \Gamma_v$.
By Lemma \ref{lem:loc diffeo} there exists a $\Gamma_v$-invariant 
neighborhood $V$ of $v$ on which $\Tilde{f}$ is a diffeomorphism
with an open subset of $\R^{n'}$.  In particular, 
$\Tilde{f}|_V$ is one to one, so

\begin{equation} \labell{as good}
 U_\gamma \cap V = U^\gamma \cap V. 
\end{equation}
\item
Suppose $\gamma \in \Gamma \ssminus \Gamma_v$.
By Lemma \ref{lem:loc diffeo}, there exist neighborhoods
$V_1$ of $v$ and $V_2$ of $\gamma \cdot v$
such that $\Tilde{f}|_{V_1}$ and $\Tilde{f}|_{V_2}$ are diffeomorphisms
with open subsets of $\R^{n'}$.

Let $V$ be a $\Gamma_v$-invariant connected neighborhood of $v$
contained in $V_1 \cap \Tilde{f}\inv(\Tilde f(V_2))$.  
For instance, we may take $V$ to be a ball centered at $v$
with respect to an invariant metric.
The map
$ u \mapsto  \left( \Tilde{f}|_V \right)\inv \Tilde{f} (\gamma \cdot u) $
from $V$ to $V$
sends each $\Gamma_v$-orbit to itself.  Applying Lemma \ref{lem:identity}
to this map, we get an element $\gamma' \in \Gamma_v$ such that
\begin{equation} \labell{as good as}
\Tilde{f}(\gamma \cdot u) = \Tilde{f}(\gamma' \cdot u)
\end{equation}
for all $u \in V$.  By the previous paragraph applied to $\gamma'$,
for each $u \in V$,
if $\Tilde{f}(\gamma' \cdot u) = \Tilde f(u)$ then $\gamma' \cdot u = u$.
By this and \eqref{as good as}, 
\begin{equation} \labell{as good as as}
 U_\gamma \cap V = U^{\gamma'} \cap V .
\end{equation}
\end{itemize}
By \eqref{as good} and \eqref{as good as as}, for every point $u \in U_\gamma$ 
there exists a neighborhood $V \subset U$ and a linear subspace
$L \subseteq \R^n$ such that
\begin{equation} \labell{intersect with L}
  U_\gamma \cap V = L \cap V .
\end{equation}
If $u \not \in \interior(U_\gamma)$
then $L$ is a proper subspace of $\R^n$,
and, by \eqref{intersect with L}, $V \cap \interior(U_\gamma) = \varnothing$.
Hence, $\interior(U_\gamma)$ and $U \ssminus \interior(U_\gamma)$
are both open in $U$.  Because $U$ is connected, either $U_\gamma = U$
or $U_\gamma$  has an empty interior.

On a sufficiently small $\Gamma$-invariant neighborhood $V$
of the origin, $\Tilde{f}$ is one to one
(by Lemma~\ref{lem:loc diffeo}), and so 
$U_\gamma \cap V = U^\gamma \cap V$ for all $\gamma \in \Gamma$.
Because $\Gamma$ acts effectively, if $\gamma$ is not the identity element,
then $U^\gamma \subset \R^n$ is a proper subspace.  Combining this
with the previous paragraph, we deduce that if $U_\gamma$ contains 
a non-empty open set then $\gamma$ is the identity.

Let $v$ be any point in $U_\gamma$.
Let $V$ be a $\Gamma_v$-invariant neighborhood of $v$ 
on which $\Tilde{f}$ is one to one.
As before, there exists $\gamma' \in \Gamma_v$ such that
$\Tilde{f}(\gamma (\gamma')\inv u) = \Tilde{f}(u)$
for all $u \in V$.
(If $\gamma \in \Gamma_v$ then we may take $\gamma'=\gamma$.)
So the set $U_{\gamma (\gamma')\inv}$ contains the open set $V$.
By the previous paragraph, $\gamma(\gamma')\inv = 1$.
So $\gamma \in \Gamma_v$.
But for $\gamma \in \Gamma_v$ we have
$V \cap U_\gamma = V \cap U^\gamma$.

We deduce that $U_\gamma = U^\gamma$.
\end{proof}

A map is called \emph{proper} if the pre-image of any compact set 
is compact.  A proper map to a first-countable Hausdorff space 
is a closed map.

\begin{proof}[Proof of (2)]
By Lemma \ref{lem:loc diffeo}, each point of $U$ has a neighborhood
on which $\Tilde{f}$ is a diffeomorphism with an open subset of $\R^{n'}$.
By this and part (1), $\Tilde{f}$ is a diffeomorphism of $U$
with an open subset of $\R^{n'}$.

Because $f(U/\Gamma) = U''/\Gamma'$ and $\Tilde{f}$ lifts $f$, 
we have $\Tilde{f}(U) \subset U''$.
Because $U$ is connected, $\Tilde{f}(U)$ is contained in a connected 
component of $U''$; denote this component by $U'$.  
Thus, $\Tilde{f}$ is an open map from $U$ to $U'$.
We need to show that its image is all of $U'$.

Consider the diagram
$$\begin{CD}
   U @> \Tilde{f} >> U' \\
   @V \pi VV @VV \pi'|_{U'} V \\
   U/\Gamma @> f >> U'' / \Gamma'.
\end{CD}$$
The projection maps $\pi \colon U \to U/\Gamma$ 
and $\pi' \colon U'' \to U''/\Gamma'$ are proper.
Because $U' \subset U''$ is closed, $\pi'|_{U'} \colon U' \to U''/\Gamma$
is also proper.  
Because the maps $\pi$, $\pi'|_{U'}$, and $f$ are proper, 
the map $\Tilde{f}$ is also proper.  
Because $\Tilde{f} \colon U \to U'$ is open and proper and $U'$ is connected, 
the image of $\Tilde{f}$ is all of $U'$.
\end{proof}

\begin{proof}[Proof of (3)]
By Part (2), $\image \Tilde{f} = U'$ is a connected component of $U''$,
and $\Tilde{f} \colon U \to U'$ is a diffeomorphism.
By Lemma \ref{lem:identity}, for each $\gamma \in \Gamma$
there exists a unique $\gamma' \in \Gamma'$ such that
\begin{equation} \labell{h}
\Tilde{f}(\gamma \cdot {\Tilde{f}}\inv(x)) = \gamma' \cdot x
\end{equation}
for all $x \in U'$.  Uniqueness implies that $\gamma \mapsto \gamma'$
defines a group homomorphism
$$ h \colon \Gamma \to \Gamma' .$$
Then \eqref{h} becomes
\begin{equation} \labell{gamma to gamma prime}
   \Tilde{f} (\gamma \cdot u) = \gamma' \cdot \Tilde{f}(u)
     \quad \text{ for all } u \in U \text{ and } \gamma' = h(\gamma),
\end{equation}
which means that $\Tilde{f}$ is $h$-equivariant.

Taking the derivatives at the origin of the left and right hand sides
of \eqref{h} gives $L \gamma = \gamma' L$,
so that the group homomorphism $h$ coincides with the map
$\gamma \mapsto L \gamma L\inv$ of Lemma~\ref{lem:loc diffeo}.
\end{proof}

We now reach the main result of this section.

\begin{Lemma}[Existence of global lifting] \labell{crucial}
Let $\Gamma \subset \GL{n}{\R}$ and $\Gamma' \subset \GL{n'}{\R}$
be finite subgroups, let $U' \subset \R^{n'}$
be a $\Gamma'$-invariant open subset, and let
$$ f \colon \R^n / \Gamma \to U' / \Gamma' $$
be a diffeomorphism.  Then there exists a smooth map
$\Tilde{f} \colon \R^n \to U'$ that lifts $f$.
Moreover, if $\Tilde{f}$ is a smooth map that lifts $f$,
then there exists a group isomorphism $h \colon \Gamma \to \Gamma'$
such that $\Tilde{f}$ is an $h$-equivariant diffeomorphism
of $\R^n$ with a connected component of $U'$. 
\end{Lemma}

\begin{proof}
For each $u \in \R^n$ and $r > 0$, let
$$ B(u,r) = \{ v \in \R^n \ | \ \| v-u \|_\Gamma < r \}$$
denote the ball around $u$ of radius $r$ with respect to the
$\Gamma$-invariant norm $\| \cdot \|_\Gamma$.
By the definition of the quotient diffeologies
there exist $\rho > 0$ and a smooth map
$$ \Tilde{f}_\rho \colon B(0,\rho) \to U' $$
that lifts $f$.  
Consider the set
\begin{multline} \labell{R}
   \calR = \{ r \ | \ \text{ there exists 
              a smooth map } \Tilde{f} \colon B(0,r) \to U' \\
\text{ that lifts $f$ and extends $\Tilde{f}_\rho$ }  . \}
\end{multline}
The set $\calR$ is non-empty, because it contains $\rho$.
Clearly, if $r \in \calR$ and $\rho\leq  r' < r$ then $r' \in \calR$.
Let us show that $\calR$ is an open subset of $[\rho,\infty)$.

Fix any $r \in \calR$.  Let 
$$ \Tilde{f} \colon B(0,r) \to U' $$
be a smooth map that lifts $f$ and extends $\Tilde{f}_\rho$.
By Lemma \ref{lifting is eq diff}, $\Tilde{f}$ is a diffeomorphism
of $B(0,r)$ with an open subset of $U'$.

Let
$$ S(0,r) = 
   \partial B(0,r) = \{ v \in \R^n \ | \ \| v \|_\Gamma = r \}.$$
By the definition of the quotient diffeologies
and by Lemma \ref{lem:loc diffeo},
for each $x \in S(0,r)$ there exist $\eps > 0$ and a lifting
$B(x,\eps) \to U'$ of $f$
which is a diffeomorphism with an open subset of $U'$.
Because $S(0,r)$ is compact, it can be covered by a finite number
of such balls, $B_i = B(x_i,\eps_i)$. 
For each $i=1,\ldots,m$, let
$$ \Tilde{f}_i \colon B_i \to U' \quad , \quad i=1,\ldots, m$$
be a lifting of $f$ which is a diffeomorphism with an open subset of $U'$.
Fix any~$i$.  By Lemma \ref{lem:identity},
there exists an element $\gamma \in \Gamma'$
such that $\Tilde{f}_i \circ \Tilde{f}\inv(x) = \gamma \cdot x$
for all $x \in \Tilde{f} ( B_i \cap B(0,r) )$.

Replacing $\Tilde{f}_i$ by $\gamma\inv \circ \Tilde{f}_i$,  we 
may assume that
\begin{equation} \labell{fi equals f}
\Tilde{f}_i = \Tilde{f} \text{ \ on \ } B_i \cap B(0,r) .
\end{equation}
Suppose that $B_i \cap B_j$ is nonempty.
By Lemma \ref{lem:identity}, there exists an element
$\gamma \in \Gamma'$ such that 
\begin{equation} \labell{eq1}
   \Tilde{f}_i \Tilde{f}_j\inv(x) = \gamma \cdot x \quad
\text{for all } x \in \Tilde{f}_j(B_i \cap B_j).
\end{equation}
Because $B_i$ and $B_j$ are centered on $S(0,r)$
and $B_i \cap B_j$ is non-empty,
the triple intersection $B(0,r) \cap B_i \cap B_j$ is also non-empty.
By \eqref{fi equals f} applied to $i$ and $j$,
\begin{equation} \labell{eq2}
\Tilde{f}_i \Tilde{f}_j\inv(x) = x
\quad \text{for all } x \in \Tilde{f}_j (B_i \cap B_j \cap B(0,r)).
\end{equation}
Because the set $\Tilde{f}_j(B_i \cap B_j \cap B(0,r))$ is open
and non-empty,
it contains an element $x$ whose $\Gamma'$-stabilizer is trivial.
For this $x$, \eqref{eq1} and \eqref{eq2} imply that $\gamma = 1$.  
So
$$ \Tilde{f} \cup \Tilde{f}_1 \cup \ldots \cup \Tilde{f}_m \ \colon
\ B(0,r) \cup B_1 \cup \ldots \cup B_m \to \Tilde{U}' $$
is a well defined lifting of $f$ which extends $\Tilde{f}_{\rho}$.
Because $B(0,r) \cup B_1 \cup \ldots \cup B_m$ contains a ball
of radius greater than $r$, this shows that the set $\calR$,
defined in \eqref{R}, is an open subset of $[\rho,\infty)$.

For each $r \in \calR$, let
$$ \Tilde{f}_r \colon B(0,r) \to \Tilde{U}'$$
be a smooth map that lifts $f$ and extends $\Tilde{f}_\rho$.
By Lemma \ref{lem:identity}, for any $r_1 < r_2$ in $\calR$ there exists
an element $\gamma \in \Gamma'$ such that
\begin{equation} \labell{one}
f_{r_2} \circ f_{r_1}\inv (x) = \gamma \cdot x
\quad \text{ for all } x \in \Tilde{f}_{r_1} (B(0,r_1)) .
\end{equation}
Because $\Tilde{f}_{r_1}$ and $\Tilde{f}_{r_2}$ both extend $\Tilde{f}_\rho$, we have
\begin{equation} \labell{two}
f_{r_2} \circ f_{r_1}\inv(x) = x
\quad \text{ for all } x \in \Tilde{f}_{r_1} (B(0,\rho)).
\end{equation}
By Lemma \ref{lem:loc diffeo}, the map $\Tilde{f}_{r_1}$ is open,
so the set $\Tilde{f}_{r_1} (B(0,\rho))$ 
contains an element $x$ whose $\Gamma'$-stabilizer is trivial.
For this $x$, \eqref{one} and \eqref{two} imply $\gamma = 1$.
So $\Tilde{f}_{r_2}\eval{B(0,r_1)} = \Tilde{f}_{r_1}.$
Then
\begin{equation} \labell{union}
   \Tilde{f} := \bigcup_{r\in\calR} \Tilde{f}_r \colon
        \ \bigcup_{r\in\calR} B(0,r) \to \Tilde{U}'
\end{equation}
is a well defined lifting of $f$ which extends $\Tilde{f}_\rho$.

This implies that if the set $\calR$ defined in \eqref{R}
is bounded then it has a maximum.  But since $\calR$ is open
in $[\rho,\infty)$, it cannot have a maximum.
Since $r \in \calR$ implies that $r' \in \calR$
for all $0 < r' < r$, we deduce that $\calR = [\rho,\infty)$.
So the domain of \eqref{union} is all of $\R^n$, and 
$$ \Tilde{f} \colon \R^n \to \Tilde{U}' $$
lifts $f$ and extends $\Tilde{f}_\rho$.

The fact that $\Tilde{f}$ is an equivariant diffeomorphism 
with a connected component of $U'$ then follows
from Lemma~\ref{lifting is eq diff}.
\end{proof}

According to the above results, a diffeomorphism
between finite linear quotients of $\R^n$ lifts
to an equivariant diffeomorphism
(by Lemma~\ref{crucial}), and this lifting
is unique up to an action of the finite linear groups
(see Lemma~\ref{lem:identity}).
The following examples that
this existence and uniqueness of equivariant liftings
may fail for differentiable functions that are not diffeomorphisms.

In both of these examples, let $\rho_n \colon \R \to \R$
denote a smooth function that takes values between $0$ and $1$,
vanishes outside the interval $[\frac{1}{n+1},\frac{1}{n}]$,
and is not always zero.

\begin{Example} \labell{example1} 
For any $\eps = (\eps_1,\eps_2,\ldots) \in \{1,-1\}^{\N}$,
let 
$$ f_{\eps} \colon \R \to \R $$
be the function
$$ f_{\eps}(x) = \begin{cases}
\eps_n e^{-\frac{1}{x}} \rho_n(x) & \text{ if }  
    \frac{1}{n+1} < x \leq \frac{1}{n} \ , \    n \in \N \\
0 & \text{ if } x > 1 \text{ or } x \leq 0 . 
\end{cases} $$
Then the functions $f_\eps$ give infinitely 
many inequivalent liftings of one function 
from $\R$ to $\R / \{ \pm 1 \}$.
Moreover, this remains true in any neighborhood of $0$.
\end{Example}

\begin{Example} \labell{example2} 
Let $f \colon \R^2 \to \R^2$ be the function
$$ f(x,y) = \begin{cases}
0 & \text{ if } r > 1 \text{ or } r = 0 \\
e^{-r} \rho_n(r) (r,0) & \text{ if } \frac{1}{n+1} < r \leq \frac{1}{n}
\text{ and $n$ is even } \\
e^{-r} \rho_n(r) (x,y) & \text{ if } \frac{1}{n+1} < r \leq \frac{1}{n}
\text{ and $n$ is odd} 
\end{cases}$$
where $r = \sqrt{x^2+y^2}$.
Fix any positive integer $m \geq 2$.
The function $f$ descends to a function $[f]$
from $\R^2/ \Z/m\Z$ to $\R^2/ \Z/m\Z$.
The function $[f]$ is differentiable because it has a smooth
lifting, namely, $f$.  On any annulus $\frac{1}{n+1} < r < \frac{1}{n}$
the function $[f]$ lifts to a smooth function
that is equivariant
with respect to a group homomorphism $h \colon \Z/m/Z \to \Z/m\Z$.
However, for even $n$ the group homomorphism $h$ must be trivial,
and for odd $n$ the group homomorphism $h$ must be the identity
homomorphism.  Consequently, there does not exist a neighborhood
of $0$ on which $[f]$ has an equivariant lifting
with respect to one group homomorphism.
\end{Example}

\section{Structure groups of diffeological orbifolds}
\labell{sec:structure groups}

Diffeology carries rich information; in particular,
we now show that the structure groups are a diffeological invariant.


\begin{Lemma} \labell{new}
Let $\Gamma \subset \GL{n}{\R}$ and $\Gamma' \subset \GL{n'}{\R}$
be finite subgroups.
Suppose that there exist open subsets $V \subset \R^n / \Gamma$
and $V' \subset \R^{n'} / \Gamma'$ that contain the origins
and a diffeomorphism $\varphi \colon V \to V'$ such that $\varphi(0) = 0$.
Then $n = n'$, and $\Gamma$ and $\Gamma'$ are conjugate in $\GL{n}{\R}$.
\end{Lemma}

\begin{proof}
Let $U_1 \subset \R^n$ and $U'_1 \subset \R^{n'}$
denote the preimages of $V$ and $V'$.

By the definitions of the quotient diffeologies
and of ``diffeomorphism",
there exists an open neighborhood $U_2$ of $0$ in $U_1$ 
and a smooth map $\Tilde{\varphi} \colon U_2 \to U_1'$ that lifts $\varphi$.
We may assume that $U_2$ is $\Gamma$-invariant
(e.g., by shrinking it to a small ball around the origin 
with respect to a $\Gamma$-invariant metric).

Let $U_2' \subset U_1'$ denote the preimage 
of $\varphi(U_2/\Gamma) \subset U_1'/\Gamma'$.
Then $\varphi$ restricts to a diffeomorphism
$ \psi \colon U_2 / \Gamma \to U_2' / \Gamma' $
such that $\psi(0) = 0$, and $\Tilde{\varphi}$ restricts to a smooth map
$ \Tilde{\psi} \colon U_2 \to U'_2 $
that lifts $\psi$.

By Lemma~\ref{lem:loc diffeo}, $n=n'$ and $\Gamma$ and $\Gamma'$
are conjugate in $\GL{n}{\R}$.
\end{proof}

\begin{Definition} \labell{structure group}
Let $X$ be a diffeological orbifold and $x \in X$ a point.
The \emph{structure group} of $X$ at $x$ is a finite subgroup
$\Gamma \subset \GL{n}{\R}$ such that there exists a
diffeomorphism from $\R^n / \Gamma$ onto a neighborhood of $x$
in $X$ that sends the origin to $x$.
By Lemma \ref{centered}, such a group exists;
by Lemma \ref{new}, such a group is unique up to conjugation
in $\GL{n}{\R}$.
\end{Definition}

%
%
%
%
%
%

A \emph{singular point} of a diffeological orbifold $X$
is a point whose structure group is non-trivial;
a \emph{regular point}, or a \emph{smooth point}, is a point 
whose structure group is trivial.
Because the structure group is unique up to conjugation,
these notions are well defined.
The set of regular points is open and dense; at these points,
$X$ is a manifold.  At a singular point, $X$ looks like
a quotient $\R^n/\Gamma$ where $\Gamma$ is a non-trivial
finite subgroup of $\GL{n}{\R}$.

To accurately state the relation of diffeological orbifolds
to Satake's V-manifolds, we need to consider orbifolds 
whose structure groups do not contain reflections.
(Also see Remark \ref{why reflection free}.)

\begin{Definition} \labell{def:reflection free}
A \emph{reflection} is a linear map $\R^n \to \R^n$ whose fixed point 
set has codimension one.
A \emph{reflection free diffeological orbifold} is a diffeological
space which is locally diffeomorphic at each point to a quotient
$\R^n / \Gamma$ for some $n$ where $\Gamma$ is a finite group
acting on $\R^n$ linearly and without reflections.
\end{Definition}

\begin{Remark} \labell{rk:reflection free}
Notice that if a subgroup of $\GL{n}{\R}$ contains no reflections 
then any subgroup of it or any group conjugate to it also contains 
no reflections.
Thus, a diffeological orbifold is reflection free if and only if 
its structure groups contain no reflections.
In particular, if $X$ is a reflection free diffeological orbifold,
$\Gamma \subset \GL{n}{\R}$ is a finite subgroup,
and $\psi \colon \R^n \to X$ induces a diffeomorphism 
from $\R^n / \Gamma$ onto an open subset of $X$, then $\Gamma$ contains
no reflections.
\end{Remark}

The following characterization will be useful to connect our notion
of a diffeological orbifold with Satake's notion of a V-manifold:

\begin{Proposition} \labell{Dorbifold II}
A diffeological space is a reflection free diffeological orbifold 
if and only if each point has a neighborhood which is diffeomorphic
to a quotient $\Tilde{U}/\Gamma$, where $\Tilde{U} \subset \R^n$ 
is a connected open subset and $\Gamma$ is a finite group 
of diffeomorphisms of $\Tilde{U}$ whose fixed point sets have 
codimension $\geq 2$.
\end{Proposition}

\begin{proof}
This follows immediately from Definition \ref{def:Dorbifold},
Example \ref{finite stabs}, and Remark \ref{rk:reflection free}.
\end{proof}

\section{V-Manifolds}
\labell{sec:Vmanifolds}

The notion of a V-manifold was introduced by Ichiro Satake
in his two papers
\emph{On a Generalization of the Notion of Manifold} \cite{satake1}
and \emph{The Gauss-Bonnet Theorem for V-Manifolds} \cite{satake2}.
Satake's definitions in \cite{satake2} slightly differ from those
in \cite{satake1}; they do lead to equivalent notions of ``V-manifold"
but this fact is not obvious.  Thus Satake himself began the tradition 
in the literature of attributing to Satake a definition which a priori 
differs from his.  We will follow the definitions of Satake's 
second paper, \cite{satake2}.

The local structure of a V-manifold is given by a
\emph{local uniformizing system} (l.u.s);
this can be thought of as a ``local chart".
The following definition is taken from \cite[p.~465--466]{satake2}.

\begin{Definition} \labell{def-lus}
Let $M$ be a Hausdorff space and $U \subset M$ an open subset.
A \emph{local uniformizing system} (l.u.s) for $U$ is a triple
$(\Tilde{U}, G, \varphi)$,
where $\Tilde{U}$ is a connected open subset of $\R^n$ for some $n$,
where $G$ is a finite group of diffeomorphisms of $\Tilde{U}$
whose fixed point sets have codimension $\geq 2$,
and where $\varphi \colon \Tilde{U} \to U$
is a map which induces a homeomorphism between $\Tilde{U}/G$ and $U$.
\end{Definition}

Local uniformizing systems are patched together by \emph{injections};
these can be thought of as the ``transition maps".
The following definition is taken from \cite[p.~466]{satake2}:

\begin{Definition} \labell{def:inject}
An \emph{injection} from an l.u.s $(\Tilde{U}, G, \varphi)$
to an l.u.s $(\Tilde{U'}, G', \varphi')$
is a diffeomorphism $\lambda$ from $\Tilde{U}$ onto an open subset
of $\Tilde{U'}$ such that
\[
	\varphi = \varphi' \circ \lambda.
\]
\end{Definition}

\begin{Remark} \labell{why reflection free}
Eventually we would like to show that diffeological orbifolds are
``the same as" V-manifolds.  However, this is not quite true.
One minor difference is that Satake only works with spaces
that are Hausdorff, whereas in our definition of a diffeological
orbifold we allow the D-topology to be non-Hausdorff.
A more serious difference is that in his definition 
of an l.u.s Satake only works with finite groups that act without reflections,
whereas we allow arbitrary finite group actions.
It appears that Satake's main usage of the reflection-free assumption
is in \cite[Lemma 1, p.~360]{satake1} and \cite[Lemma 1, 
p.~466]{satake2}, which assert this:
\begin{Lemma*} 
Let \function{\lambda,\mu}{\Tilde{U}}{\Tilde{U'}} be two injections.
Then there exists a unique $\sigma'\in G'$
such that $\mu = \sigma' \circ \lambda$.
\end{Lemma*}
However, this lemma remains true even if the finite groups contain
reflections.
This more general version of the lemma appeared in \cite[Appendix]{MP}.
It can also be proved in a way similar to our proof of Lemma
\ref{lem:identity}. 
\end{Remark}

\begin{figure}[ht]
\setlength{\unitlength}{.00001in}
\begin{center}
\setlength{\unitlength}{0.0005in}
\begingroup\makeatletter\ifx\SetFigFont\undefined%
\gdef\SetFigFont#1#2#3#4#5{%
  \reset@font\fontsize{#1}{#2pt}%
  \fontfamily{#3}\fontseries{#4}\fontshape{#5}%
  \selectfont}%
\fi\endgroup%
{\renewcommand{\dashlinestretch}{30}
\begin{picture}(5478,6334)(0,-10)
\put(2261,1512){\ellipse{2100}{2100}}
\put(3161,1512){\ellipse{2100}{2100}}
\put(2711,1512){\ellipse{618}{618}}
\put(2711,1512){\whiten\ellipse{80}{80}}
\put(2711,1512){\ellipse{80}{80}}
\put(2711,6012){\ellipse{600}{600}}
\put(4361,4512){\ellipse{2100}{2100}}
\put(1061,4512){\ellipse{2106}{2106}}
\put(3836,4512){\ellipse{600}{600}}
\put(1511,4512){\ellipse{600}{600}}
\path(2711,5702)(2711,1852)
\path(2681.000,1972.000)(2711.000,1852.000)(2741.000,1972.000)
\path(4061,3462)(3706,2432)
\path(3718.108,2555.099)(3706.000,2432.000)(3774.613,2534.919)
\path(1361,3462)(1716,2432)
\path(1647.387,2534.919)(1716.000,2432.000)(1703.892,2555.099)
\path(2411,5862)(1736,4737)
\path(1772.015,4855.334)(1736.000,4737.000)(1823.464,4824.464)
\path(2996,5862)(3611,4737)
\path(3523.536,4824.464)(3611.000,4737.000)(3574.985,4855.334)
\path(161,2862)(5261,2862)(5261,12)
	(161,12)(161,2862)
\put(2541,1512){\makebox(0,0)[lb]{\smash{{{\SetFigFont{12}{14.4}{\rmdefault}{\mddefault}{\updefault}$\sss p$}}}}}
\put(951,1512){\makebox(0,0)[lb]{\smash{{{\SetFigFont{12}{14.4}{\rmdefault}{\mddefault}{\updefault}$\sss U_1$}}}}}
\put(4286,1512){\makebox(0,0)[lb]{\smash{{{\SetFigFont{12}{14.4}{\rmdefault}{\mddefault}{\updefault}$\sss U_2$}}}}}
\put(2806,1077){\makebox(0,0)[lb]{\smash{{{\SetFigFont{12}{14.4}{\rmdefault}{\mddefault}{\updefault}$\sss U_3$}}}}}
\put(5166,5187){\makebox(0,0)[lb]{\smash{{{\SetFigFont{12}{14.4}{\rmdefault}{\mddefault}{\updefault}$\sss \Tilde{U}_2$}}}}}
\put(1,5187){\makebox(0,0)[lb]{\smash{{{\SetFigFont{12}{14.4}{\rmdefault}{\mddefault}{\updefault}$\sss \Tilde{U}_1$}}}}}
\put(3086,6012){\makebox(0,0)[lb]{\smash{{{\SetFigFont{12}{14.4}{\rmdefault}{\mddefault}{\updefault}$\sss \Tilde{U}_3$}}}}}
\put(3276,5412){\makebox(0,0)[lb]{\smash{{{\SetFigFont{12}{14.4}{\rmdefault}{\mddefault}{\updefault}$\sss \lambda_2$}}}}}
\put(1896,5412){\makebox(0,0)[lb]{\smash{{{\SetFigFont{12}{14.4}{\rmdefault}{\mddefault}{\updefault}$\sss \lambda_1$}}}}}
\put(2786,3462){\makebox(0,0)[lb]{\smash{{{\SetFigFont{12}{14.4}{\rmdefault}{\mddefault}{\updefault}$\sss \varphi_3$}}}}}
\put(4961,2637){\makebox(0,0)[lb]{\smash{{{\SetFigFont{12}{14.4}{\rmdefault}{\mddefault}{\updefault}$\sss V$}}}}}
\put(1211,3087){\makebox(0,0)[lb]{\smash{{{\SetFigFont{12}{14.4}{\rmdefault}{\mddefault}{\updefault}$\sss \varphi_1$}}}}}
\put(3986,3087){\makebox(0,0)[lb]{\smash{{{\SetFigFont{12}{14.4}{\rmdefault}{\mddefault}{\updefault}$\sss \varphi_2$}}}}}
\end{picture}
}
\end{center}
\caption{Defining family \`{a} la Satake}
\label{fig:injections}
\end{figure}
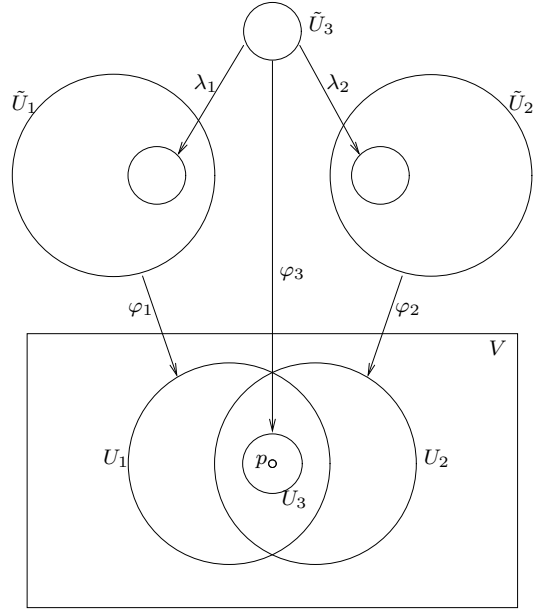

A V-manifold structure is given by a \emph{defining family};
this can be thought of as an ``atlas". The following two definitions
are taken from \cite[p.~467]{satake2}.

\begin{Definition} \labell{def:defining family}
Let $M$ be a Hausdorff space.  A \emph{defining family} on $M$
is a family $\calF$ of l.u.s's for open subsets of $M$,
satisfying conditions (1) and (2) below.
An open subset $U \subset M$ is said to be $\calF$-\emph{uniformized}
if there exists an l.u.s.\ $(\Tilde{U},G,\varphi)$ in $\calF$
such that $\varphi(\Tilde{U}) = U$.
\begin{enumerate}
\item
Every point in $M$ is contained in at least one $\calF$-uniformized open
set.  If a point $p$ is contained in two $\calF$-uniformized open sets
$U_1$ and $U_2$ then there exists an $\calF$-uniformized open set $U_3$
such that $p \in U_3 \subset U_1 \cap U_2$.
\item
If $(\Tilde{U},G,\varphi)$ and $(\Tilde{U'},G,\varphi')$
are l.u.s's in $\calF$
and $\varphi(\Tilde{U}) \subset \varphi'(\Tilde{U'})$,
then there exists an injection
\function{\lambda}{\Tilde{U}}{\Tilde{U'}}.
\end{enumerate}
\end{Definition}

See Figure \ref{fig:injections}.

One might naively attempt to declare two defining families
to be equivalent if their union is a defining family.  However, 
this is not an equivalence relation:

\begin{Example} \labell{bad union}
Consider the annulus
$$ M = \{ (x,y) \in \R^2 \mid 0 < x^2 + y^2 < 1 \}.$$
Consider the following two l.u.s's.
The first consists of the annulus with the trivial group action
and the identity map.  The second consists of the annulus
with the two element group acting by $(x,y) \mapsto \pm (x,y)$
and with the map
$\varphi (r\cos\theta , r\sin\theta) = (r\cos2\theta , r\sin2\theta)$.

Let $\calF$ consist of the first l.u.s, and let $\calF'$ consist
of the second l.u.s.  Let $\calF''$ consist of all l.u.s's
consisting of a disk that is contained in the annulus, with
the trivial group action and the identity map.
Then $\calF$, $\calF'$, and $\calF''$ are defining families
on $M$. The unions $\calF \cup \calF''$ and $\calF' \cup \calF''$
are defining families, but there is no defining family that contains
both $\calF$ and $\calF'$.
\end{Example}

\begin{Definition} \labell{def:equivalent F}
Two defining families, $\calF$ and $\calF'$,
are \emph{directly equivalent} if there exists a third defining family
that contains both $\calF$ and $\calF'$.
Two defining families $\calF$ and $\calF'$ are \emph{equivalent}
if there exists a chain of direct equivalences
starting with $\calF$ and ending with $\calF'$.
\end{Definition}

The following definition is taken from
\cite[p.~467, Definition 1 and footnote 1]{satake2}:

\begin{Definition} \labell{def:Vmanifold}
A \emph{V-manifold} is a Hausdorff space $M$
equipped with an equivalence class of defining families.
\end{Definition}

In what follows, whenever we choose a defining family
on a V-manifold, we assume that this family belongs to the
given equivalence class, unless otherwise stated.

\section{V-manifolds as diffeological orbifolds}
\labell{sec:Vmanifolds as Dorbifolds}

In this section we describe the natural correspondence
between Satake's V-manifolds and diffeological orbifolds.

\begin{proposition} \labell{Vmanifold is Dorbifold}
Let $M$ be a Hausdorff topological space.
A defining family $\calF$ on~$M$ determines a diffeology,
namely, the diffeology generated by the maps 
$\function{\varphi}{\Tilde{U}}{U}$,
for all $(\Tilde{U},G,\varphi) \in \calF$.
\begin{enumerate}
\item
Let $\calF$ be a defining family on $M$.
Then, equipped with the diffeology generated by $\calF$,
$M$ is a reflection-free diffeological orbifold.
\item
Equivalent defining families $\calF$, $\calF'$ on $M$
generate the same diffeology.
\end{enumerate}
\end{proposition}

The \emph{natural diffeology} on a V-manifold
is the one determined by any defining family;
see Definition~\ref{def:Vmanifold} 
and Proposition~\ref{Vmanifold is Dorbifold}.

\begin{Theorem} \labell{main theorem}
A V-manifold, with its natural diffeology, becomes
a Hausdorff reflection-free diffeological orbifold.
The diffeology uniquely determines the V-manifold structure.
Every Hausdorff reflection free diffeological orbifold arises in this way.
\end{Theorem}

The rest of this section is devoted to the proof of 
Proposition~\ref{Vmanifold is Dorbifold}
and Theorem~\ref{main theorem}.

\begin{Lemma} \labell{inclusion diffeo}
Let $M$ be a Hausdorff topological space, let $U \subset U' \subset M$
be open subsets, let $(\Tilde{U},G,\varphi)$ and $(\Tilde{U}',G',\varphi')$
be l.u.s's for $U$ and $U'$, 
and let $\lambda \colon \Tilde{U} \to \Tilde{U}'$ be an injection. Then
\begin{enumerate}
\item
The image of $\lambda$ is a connected component of $(\varphi')\inv(U)$.
\item
$\lambda$ descends to a diffeomorphism from $\Tilde{U}/G$ to an open subset
of $\Tilde{U}'/G'$.
\end{enumerate}
\end{Lemma}

\begin{proof}
%
%
From the definition of an l.u.s and of an injection,
it follows that $\lambda$ descends
to a map $\ol{\lambda} \colon \tilde{U}/G \to \Tilde{U}'/G'$
which is a homeomorphism with the open subset $(\varphi')\inv(U) / G'$
of $\Tilde{U}'/G'$.

Now consider the commuting diagram
$$ \begin{CD}
\Tilde{U} @> \lambda >> (\varphi')\inv(U) \\
@VVV @VVV \\
\Tilde{U}/G @> \ol{\lambda} >> (\varphi')\inv(U)/G' .
\end{CD} $$
Because the quotient maps and the homeomorphism $\lambda$ are proper,
$\lambda$ is proper, so its image is closed.  But, by assumption,
the image of $\lambda$ is also open and connected, so this image
must be a connected component of $(\varphi')\inv(U)$.

Finally, consider the diagram
$$ \begin{CD}
\Tilde{U} @> \lambda >> \Tilde{U}' \\
@VVV @VVV \\
\Tilde{U}/G @> \ol{\lambda} >> \Tilde{U}'/G' .
\end{CD} $$
Because $\lambda$ is an open inclusion and the vertical arrows
generate the quotient diffeologies, the map $\ol{\lambda}$
is a diffeomorphism with its image.
\end{proof}


\begin{Corollary} \labell{local diffeos}
Let $M$ be a Hausdorff space and $\calF$ a defining family;
equip $M$ with the diffeology generated by $\calF$.
Let $(\Tilde{U},G,\varphi)$ be an l.u.s in $\calF$.
Then $U := \varphi(\Tilde{U}) \subset M$ is D-open,
and the homeomorphism $\function{\ol{\varphi}}{\Tilde{U}/G}{U}$
induced by $\varphi$ is a diffeomorphism of diffeological spaces.
\end{Corollary}


Part (1) of Proposition~\ref{Vmanifold is Dorbifold}
follows from Proposition~\ref{Dorbifold II}
and Corollary~\ref{local diffeos}.
Part (2) follows from Lemma~\ref{D cover}, 
Lemma~\ref{inclusion diffeo}, and Corollary~\ref{local diffeos}.

\begin{Lemma} \labell{family of balls}
Let $M$ be a Hausdorff reflection-free diffeological orbifold.
Let $\calU$ be an open cover of $M$. Let $\calF$ consist of all
the triples $(\R^n,G,\psi)$ where $G \subset \GL{n}{\R}$ is a finite
subgroup, $\psi \colon \R^n \to M$ induces a diffeomorphism from
$\R^n/G$
onto an open subset of $M$, and there exists $U \in \calU$ such that
$\psi(\R^n) \subset U$.  Then $\calF$ is a defining family,
and the set of $\calF$-uniformized open sets is a basis to the
topology of $M$.
\end{Lemma}

\begin{proof}
By Remark \ref{rk:reflection free}, for each $(\R^n,G,\psi) \in \calF$ 
the group $G$ is reflection free.
Thus, $(\R^n,G,\psi)$ is an l.u.s.
By Remark \ref{part of orbifold} and Lemma \ref{centered},
the $\calF$-uniformized open sets form a basis to the topology of $M$.
So $\calF$ satisfies Condition (1) of Definition
\ref{def:defining family} of a defining family.
Condition (2) follows
from Lemma \ref{crucial}. See Figure \ref{fig:globallift}.
\end{proof}

\begin{figure}[ht]
\begin{center}
\setlength{\unitlength}{0.0005in}
\begingroup\makeatletter\ifx\SetFigFont\undefined%
\gdef\SetFigFont#1#2#3#4#5{%
  \reset@font\fontsize{#1}{#2pt}%
  \fontfamily{#3}\fontseries{#4}\fontshape{#5}%
  \selectfont}%
\fi\endgroup%
{\renewcommand{\dashlinestretch}{30}
\begin{picture}(4975,3710)(0,-10)
\path(300,3312)(1200,3312)(1200,2712)
	(300,2712)(300,3312)
\path(2325,762)(2100,1062)(2625,1137)
	(2850,837)(2325,762)
\path(825,2712)(2325,1137)
\path(2220.517,1203.207)(2325.000,1137.000)(2263.966,1244.586)
\path(3950,2412)(3150,1212)
\path(3191.603,1328.487)(3150.000,1212.000)(3241.526,1295.205)
\dottedline{45}(3750,3312)(4650,3312)(4650,2712)
	(3750,2712)(3750,3312)
\dottedline{45}(1200,3012)(3750,3012)
\path(3630.000,2982.000)(3750.000,3012.000)(3630.000,3042.000)
\put(4200,3012){\ellipse{1350}{1350}}
\path(900,1587)(827,1560)(760,1532)
	(697,1502)(639,1472)(586,1441)
	(537,1410)(493,1379)(453,1349)
	(416,1319)(383,1289)(353,1260)
	(327,1232)(303,1205)(281,1179)
	(262,1153)(245,1128)(229,1104)
	(216,1081)(204,1058)(193,1036)
	(183,1015)(175,994)(167,973)
	(161,953)(155,932)(150,912)
	(146,892)(143,871)(140,851)
	(138,830)(137,809)(137,788)
	(139,766)(141,743)(144,720)
	(149,696)(155,671)(163,645)
	(172,619)(184,592)(199,564)
	(216,535)(236,505)(259,475)
	(285,445)(316,414)(351,383)
	(391,352)(435,322)(484,292)
	(539,264)(600,237)(657,215)
	(718,195)(781,176)(847,159)
	(914,143)(982,128)(1050,115)
	(1119,103)(1187,92)(1255,82)
	(1321,73)(1387,64)(1451,57)
	(1514,51)(1576,45)(1636,40)
	(1694,35)(1751,31)(1807,28)
	(1861,25)(1914,22)(1965,20)
	(2016,18)(2065,16)(2114,15)
	(2162,14)(2210,13)(2257,12)
	(2303,12)(2350,12)(2397,12)
	(2444,12)(2491,13)(2539,14)
	(2587,15)(2637,16)(2687,18)
	(2738,20)(2791,22)(2845,25)
	(2900,28)(2958,31)(3016,35)
	(3077,40)(3139,45)(3203,51)
	(3269,57)(3336,64)(3406,73)
	(3476,82)(3548,92)(3621,103)
	(3695,115)(3770,128)(3844,143)
	(3919,159)(3992,176)(4064,195)
	(4133,215)(4200,237)(4273,264)
	(4340,292)(4403,322)(4461,352)
	(4514,383)(4563,414)(4607,445)
	(4647,475)(4684,505)(4717,535)
	(4747,564)(4773,592)(4797,619)
	(4819,645)(4838,671)(4855,696)
	(4871,720)(4884,743)(4896,766)
	(4907,788)(4917,809)(4925,830)
	(4933,851)(4939,871)(4945,892)
	(4950,912)(4954,932)(4957,953)
	(4960,973)(4962,994)(4963,1015)
	(4963,1036)(4961,1058)(4959,1081)
	(4956,1104)(4951,1128)(4945,1153)
	(4937,1179)(4928,1205)(4916,1232)
	(4901,1260)(4884,1289)(4864,1319)
	(4841,1349)(4815,1379)(4784,1410)
	(4749,1441)(4709,1472)(4665,1502)
	(4616,1532)(4561,1560)(4500,1587)
	(4443,1609)(4382,1629)(4319,1648)
	(4253,1665)(4186,1681)(4118,1696)
	(4050,1709)(3981,1721)(3913,1732)
	(3845,1742)(3779,1751)(3713,1760)
	(3649,1767)(3586,1773)(3524,1779)
	(3464,1784)(3406,1789)(3349,1793)
	(3293,1796)(3239,1799)(3186,1802)
	(3135,1804)(3084,1806)(3035,1808)
	(2986,1809)(2938,1810)(2890,1811)
	(2843,1812)(2797,1812)(2750,1812)
	(2703,1812)(2656,1812)(2609,1811)
	(2561,1810)(2513,1809)(2463,1808)
	(2413,1806)(2362,1804)(2309,1802)
	(2255,1799)(2200,1796)(2142,1793)
	(2084,1789)(2023,1784)(1961,1779)
	(1897,1773)(1831,1767)(1764,1760)
	(1694,1751)(1624,1742)(1552,1732)
	(1479,1721)(1405,1709)(1330,1696)
	(1256,1681)(1181,1665)(1108,1648)
	(1036,1629)(967,1609)(900,1587)
\put(0,312){\makebox(0,0)[lb]{\smash{{{\SetFigFont{12}{14.4}{\rmdefault}{\mddefault}{\updefault}$\sss X$}}}}}
\path(3225,687)(3198,664)(3167,642)
	(3134,623)(3099,605)(3063,590)
	(3026,576)(2991,564)(2956,553)
	(2922,544)(2890,536)(2858,529)
	(2828,522)(2799,517)(2770,512)
	(2742,507)(2715,503)(2687,499)
	(2660,496)(2632,493)(2603,491)
	(2573,488)(2542,486)(2509,485)
	(2474,484)(2437,483)(2398,484)
	(2357,485)(2313,487)(2267,491)
	(2220,496)(2171,503)(2122,512)
	(2073,523)(2025,537)(1977,555)
	(1934,575)(1894,597)(1857,619)
	(1825,641)(1796,664)(1771,686)
	(1748,707)(1728,727)(1711,746)
	(1695,765)(1681,783)(1669,800)
	(1657,817)(1647,833)(1638,849)
	(1629,866)(1621,883)(1613,900)
	(1606,918)(1599,937)(1594,957)
	(1589,977)(1585,1000)(1583,1023)
	(1583,1048)(1585,1074)(1589,1102)
	(1598,1129)(1610,1158)(1627,1185)
	(1650,1212)(1677,1235)(1708,1257)
	(1741,1276)(1776,1294)(1812,1309)
	(1849,1323)(1884,1335)(1919,1346)
	(1953,1355)(1985,1363)(2017,1370)
	(2047,1377)(2076,1382)(2105,1387)
	(2133,1392)(2160,1396)(2188,1400)
	(2215,1403)(2243,1406)(2272,1408)
	(2302,1411)(2333,1413)(2366,1414)
	(2401,1415)(2438,1416)(2477,1415)
	(2518,1414)(2562,1412)(2608,1408)
	(2655,1403)(2704,1396)(2753,1387)
	(2802,1376)(2850,1362)(2895,1345)
	(2936,1326)(2975,1306)(3009,1285)
	(3041,1264)(3069,1243)(3094,1222)
	(3117,1202)(3137,1183)(3154,1164)
	(3170,1146)(3184,1129)(3197,1112)
	(3208,1096)(3219,1080)(3229,1065)
	(3238,1049)(3246,1034)(3254,1018)
	(3261,1002)(3268,985)(3274,968)
	(3279,950)(3284,930)(3288,910)
	(3291,888)(3292,866)(3292,842)
	(3289,817)(3284,791)(3276,765)
	(3263,738)(3247,712)(3225,687)
\put(0,2937){\makebox(0,0)[lb]{\smash{{{\SetFigFont{12}{14.4}{\rmdefault}{\mddefault}{\updefault}$\sss \R^n$}}}}}
\put(2250,3062){\makebox(0,0)[lb]{\smash{{{\SetFigFont{12}{14.4}{\rmdefault}{\mddefault}{\updefault}$\sss \lambda$}}}}}
\put(4850,3312){\makebox(0,0)[lb]{\smash{{{\SetFigFont{12}{14.4}{\rmdefault}{\mddefault}{\updefault}$\sss U$}}}}}
\put(1100,2112){\makebox(0,0)[lb]{\smash{{{\SetFigFont{12}{14.4}{\rmdefault}{\mddefault}{\updefault}$\sss \varphi$}}}}}
\put(3825,2112){\makebox(0,0)[lb]{\smash{{{\SetFigFont{12}{14.4}{\rmdefault}{\mddefault}{\updefault}$\sss \varphi'$}}}}}
\end{picture}
}
\caption{} 
\label{fig:globallift}
\end{center}
\end{figure}

\begin{proposition} \labell{families are equivalent}
Let $M$ be a Hausdorff topological space.
If two defining families on $M$ generate the same diffeology
then they are equivalent.
\end{proposition}

\begin{proof}
Let $\calF$ and $\calF'$ be defining families on $M$
that generate the same diffeology.
Let $\calF''$ be the defining family that is constructed
as in Lemma \ref{family of balls} with respect to the open cover
$\calU = \{ U \cap U' \ | $
$U$ is $\calF$-uniformized and $U'$ is $\calF'$-uniformized $ \}$.
Let us show that $\calF \cup \calF''$ is a defining family.

Because the $\calF''$-uniformized sets form a basis to the topology,
$\calF \cup \calF''$ satisfies Condition (1)
of Definition \ref{def:defining family}.  To prove Condition (2),
suppose that $(\Tilde{U}_i,G_i,\varphi_i)$, for $i=1,2$,
are l.u.s's in $\calF \cup \calF''$ 
with $\varphi(\Tilde{U}_1) \subset \varphi(\Tilde{U}_2)$.
We will now show that there exists an injection from
$(\Tilde{U}_1,G_1,\varphi_1)$ to $(\Tilde{U}_2,G_2,\varphi_2)$.
If both l.u.s's are in $\calF$ then the fact that $\calF$ is a defining
family implies that there exists such an injection.
If $(\Tilde{U}_1,G_1,\varphi_1) \in \calF''$ then Lemma \ref{crucial}
implies that there exists such an injection. Suppose now that
$(\Tilde{U}_1,G_1,\varphi_1) \in \calF$ and
$(\Tilde{U}_2,G_2,\varphi_2) \in \calF''$.

By the definition of $\calF''$ there exists an l.u.s
$(\Tilde{U}_3,G_3,\varphi_3) \in \calF$ such that 
$\varphi(\Tilde{U}_2) \subset \varphi(\Tilde{U}_3)$.
Because $\calF$ is a defining family, there exists an open embedding
$\alpha \colon \Tilde{U}_1 \to \Tilde{U}_3$ such that
$\varphi_3 \circ \alpha = \varphi_1$.
By Lemma \ref{crucial} there exists a diffeomorphism
$\beta$ from $\Tilde{U}_2$ to a connected component of
$\varphi_3\inv(\varphi_2(\Tilde{U}_2))$ such that
$\varphi_3 \circ \beta = \varphi_2$.  
By possibly composing with an element of $G_3$ we may assume
that this connected component is the one that contains
$\alpha(\Tilde{U}_1)$.  Then
$\lambda := \beta\inv \circ \alpha \colon \Tilde{U}_1 \to \Tilde{U}_2$
is an injection, as required.
See Figure \ref{fig:proof}.

\begin{figure}[ht]
\setlength{\unitlength}{0.0005in}
\begingroup\makeatletter\ifx\SetFigFont\undefined%
\gdef\SetFigFont#1#2#3#4#5{%
  \reset@font\fontsize{#1}{#2pt}%
  \fontfamily{#3}\fontseries{#4}\fontshape{#5}%
  \selectfont}%
\fi\endgroup%
{\renewcommand{\dashlinestretch}{30}
\begin{picture}(7591,5416)(0,-10)
\put(3233,907){\ellipse{6450}{1800}}
\put(2708,907){\ellipse{600}{450}}
\put(533,4207){\ellipse{670}{670}}
\put(5933,4207){\ellipse{3300}{1950}}
\put(5858,4282){\ellipse{750}{450}}
\path(2033,1357)(1733,457)(4133,457)
	(3908,1357)(2033,1357)
\path(1883,3757)(3683,3757)(3683,4807)
	(1883,4807)(1883,3757)
\path(5033,4807)(6833,4807)(6833,3757)
	(5033,3757)(5033,4807)
\path(533,3907)(2483,1057)
\path(2390.479,1139.096)(2483.000,1057.000)(2439.997,1172.977)
\path(2933,3757)(2933,1357)
\path(2903.000,1477.000)(2933.000,1357.000)(2963.000,1477.000)
\path(5933,3232)(5633,1507)
\path(5624.005,1630.366)(5633.000,1507.000)(5683.117,1620.085)
\put(3114.739,1700.767){\arc{7329.163}{4.0138}{5.4422}}
\path(5447.698,4487.977)(5558.000,4432.000)(5487.017,4533.298)
\path(5633,3757)(3683,1357)
\path(3735.388,1469.051)(3683.000,1357.000)(3781.955,1431.216)
\put(4208,5257){\makebox(0,0)[lb]{\smash{{{\SetFigFont{12}{14.4}{\rmdefault}{\mddefault}{\updefault}$\sss \alpha$}}}}}
\path(3683,4207)(5033,4207)
\path(4913.000,4177.000)(5033.000,4207.000)(4913.000,4237.000)
\put(383,4602){\makebox(0,0)[lb]{\smash{{{\SetFigFont{12}{14.4}{\rmdefault}{\mddefault}{\updefault}$\sss \Tilde{U}_1$}}}}}
\put(1208,2522){\makebox(0,0)[lb]{\smash{{{\SetFigFont{12}{14.4}{\rmdefault}{\mddefault}{\updefault}$\sss \varphi_1$}}}}}
\put(2483,4877){\makebox(0,0)[lb]{\smash{{{\SetFigFont{12}{14.4}{\rmdefault}{\mddefault}{\updefault}$\sss \Tilde{U}_2=\R^n$}}}}}
\put(2683,2597){\makebox(0,0)[lb]{\smash{{{\SetFigFont{12}{14.4}{\rmdefault}{\mddefault}{\updefault}$\sss \varphi_2$}}}}}
\put(7163,4882){\makebox(0,0)[lb]{\smash{{{\SetFigFont{12}{14.4}{\rmdefault}{\mddefault}{\updefault}$\sss \Tilde{U}_3$}}}}}
\put(6908,4207){\makebox(0,0)[lb]{\smash{{{\SetFigFont{12}{14.4}{\rmdefault}{\mddefault}{\updefault}$\sss \Tilde{U}_4$}}}}}
\put(5858,2257){\makebox(0,0)[lb]{\smash{{{\SetFigFont{12}{14.4}{\rmdefault}{\mddefault}{\updefault}$\sss \varphi_3$}}}}}
\put(4568,2257){\makebox(0,0)[lb]{\smash{{{\SetFigFont{12}{14.4}{\rmdefault}{\mddefault}{\updefault}$\sss \varphi_4$}}}}}
\put(3983,4282){\makebox(0,0)[lb]{\smash{{{\SetFigFont{12}{14.4}{\rmdefault}{\mddefault}{\updefault}$\sss \beta$}}}}}
\end{picture}
}
\caption{}
\labell{fig:proof}
\end{figure}

We have shown that $\calF \cup \calF''$ is a defining family.
Applying a similar argument to $\calF'$ and $\calF''$,
we get a chain of direct equivalences: $\calF$, $\calF''$, $\calF'$.  
By Definition \ref{def:equivalent F}, 
the defining families $\calF$ and $\calF'$ are equivalent.
\end{proof}

\begin{proposition} \labell{Dorbifold is Vmanifold}
Let $M$ be a Hausdorff reflection free diffeological orbifold.
Then there exists a defining family $\calF$ on $M$
which generates its diffeology.
\end{proposition}

\begin{proof}
Let $\calF$ consist of all the triples $(\R^n,G,\varphi)$
where $G \subset \GL{n}{\R}$ is a finite subgroup
and where $\varphi \colon \R^n \to M$ induces a diffeomorphism
from $\R^n/G$, equipped with its quotient diffeology,
to an open subset of $M$, equipped with its subset diffeology.
By Lemma \ref{family of balls} (with any $\calU$),
$\calF$ is a defining family.
By Corollary \ref{local diffeos} and the definition of $\calF$,
the family $\calF$ generates the diffeology on $M$.
\end{proof}

Theorem~\ref{main theorem} follows immediately from 
Propositions~\ref{families are equivalent}
and~\ref{Dorbifold is Vmanifold}.

\section{Locality of V-manifold structure}
\labell{sec:locality}

Let $X$ be a topological space and $Y \subset X$ an open subset.
Given a defining family $\calF$ on $X$, its restriction
$\calF|_Y$ to $Y$ is defined to be the set of all l.u.s's
of the form $(\Tilde{U}_Y,G,\varphi|_{\Tilde{U}_Y})$
where $(\Tilde{U},G,\varphi)$ is an l.u.s on $X$
and where $\Tilde{U}_Y$ is a connected component 
of $\Tilde{U} \cap \varphi\inv(Y)$.
In this way, a V-manifold structure on $X$ naturally restricts
to a V-manifold structure on $Y$.

Given two defining families $\calF_1$ and $\calF_2$, if the restrictions
$\calF_1|_U$ and $\calF_2|_U$ are equivalent for elements
$U$ of an open covering $\calU$ of $M$, a priori it is not clear
that $\calF_1$ and $\calF_2$ themselves are equivalent.
Also, given defining families $\calF_U$ for $U \in \calU$,
if the restrictions $\calF_U|_{U \cap U'}$ 
and $\calF_{U'}|_{U \cap U'}$ are equivalent 
for every $U,U' \in \calU$, a priori it is not clear
how to patch together the families $\calF_U$ into a defining
family on $M$.  However, these locality properties
do follow, immediately, from the fact that a V-manifold
structure is determined by the diffeology that it generates,
yielding the following result:


\begin{Proposition} \labell{locality}
Let $M$ be a Hausdorff topological space and $\calU$ an open cover.

A V-manifold structure on $M$ is uniquely determined by its restrictions
to the open sets $U$ in $\calU$.

Suppose that we are given for each $U \in \calU$ a V-manifold structure
on $U$.  Suppose that, for each $U$ and $U'$ in $\calU$, the restrictions
to $U \cap U'$ of the V-manifold structures on $U$ and on $U'$ coincide.  
Then there exists a unique V-manifold structure on $M$ whose restriction
to each $U \in \calU$ is the given V-manifold structure on~$U$.
\end{Proposition}

\section{Relation with Haefliger's definition}
\labell{sec:Haefliger}

In section 4 of \cite{haefliger}, Haefliger gives this definition.

\begin{quotation}
The structure of a V-manifold on $Q$ is given by an open covering $U_i$,
each $U_i$ being the quotient of a ball $\tU_i$ in $\R^q$
by a finite group $G_i$ of diffeomorphisms.
The projection $\pi_i \colon \tU_i \to U_i = \tU_i / G_i$
is called a uniformizing map.  If $\tu_i$ and $\tu_j$ are two points
in $\tU_i$ and $\tU_j$ such that $\pi_i(\tu_i) = \pi_j(\tu_j)$,
one assumes that there exists a diffeomorphism $g_{ji}$
from a neighborhood of $\tu_i$ onto a neighborhood of $\tu_j$
such that $\pi_j \circ g_{ji} = \pi_i$ on a neighborhood of $\tu_i$.
\end{quotation}

Thus, Haefliger works with triples $(\tU_i,G_i,\pi_i)$,
where $G_i$ is a finite group of diffeomorphisms of $\tU_i$
and the map $\pi_i \colon \tU_i \to Q$ induces a homeomorphism 
$\pibar_i \colon \tU_i/G_i \to U_i$ onto an open subset $U_i$ of $Q$.
He gives a compatibility condition for such triples.

\begin{Theorem}
Let $\{ ( \tU_i , G_i , \pi_i ) \}$ be a collection
that satisfies Haefliger's compatibility condition
such that the sets $U_i$ cover $Q$.
Then $Q$, equipped with the diffeology generated
by the maps $\pi_i$, is a diffeological orbifold.
Two such collections generate the same diffeology on $Q$
if and only if their union satisfies Haefliger's
compatibility condition.  Every diffeological orbifold
is obtained in this way. 
\end{Theorem}

\begin{proof}
Haefliger's compatibility condition implies
that the transition functions
${\pibar_j}\inv \circ \pibar_i$ are differentiable
in the diffeological sense.  By Lemma~\ref{local liftings}
and the slice theorem, if ${\pibar_j}\inv \circ \pibar_i$
is a diffeomorphism in the diffeological sense,
then $(\tU_i,G_i,\pi_i)$ and $(\tU_j,G_j,\pi_j)$
are compatible in Haefliger's sense.  These two facts
imply the theorem.
\end{proof}

Thus, an orbifold structure in Haefliger's sense
is the same thing as a diffeological orbifold;
hence, in the reflection-free case, it is also
the same thing as Satake's V-manifold.
A priori this is not obvious from the definitions
in~\cite{haefliger} and~\cite{satake1,satake2}.

\end{document}